\documentstyle[11pt,amsfonts,fullpage]{amsart}

\def\endproof{\qed \medskip}
\def\blacksquare{\hbox to .60em{\vrule width .60em height .60em}}

\newtheorem{theorem}{Theorem}[section]

\newtheorem{corollary}[theorem]{Corollary}

\newtheorem{lemma}[theorem]{Lemma}

\newtheorem{proposition}[theorem]{Proposition}

\newtheorem{remark}[theorem]{Remark}

\begin{document}

\title[]{Orbifold Compactness for Spaces of Riemannian Metrics and 
Applications}

\author[]{Michael T. Anderson}

\thanks{Partially supported by NSF Grant DMS 0305865}

\maketitle

\setcounter{section}{0}

\section{Introduction}
\setcounter{equation}{0}

 The Cheeger-Gromov compactness theorem, cf. [C], [G], [CGv] states that 
the space of Riemannian $n$-manifolds $(M^{n}, g)$ satisfying the bounds
\begin{equation} \label{e1.1}
|R| \leq  \Lambda , \ vol \geq  v_{0}, \ diam \leq  D, 
\end{equation}
is precompact in the $C^{1,\alpha}$ topology. Here $R$ denotes the 
Riemann curvature tensor, $vol$ the volume and $diam$ the diameter of $(M, g)$. 
Thus, for any sequence of metrics $g_{i}$ on $n$-manifolds $M_{i}$ 
satisfying (1.1), there is a subsequence, also called $g_{i},$ and 
diffeomorphisms $\phi_{i}: M_{\infty} \rightarrow  M_{i}$ such that the 
metrics $\phi_{i}^{*}g_{i}$ converge in the $C^{1,\alpha}$ topology to 
a limit metric $g_{\infty}$ on $M_{\infty},$ for any $\alpha  < $ 1. In 
particular, there are only finitely many diffeomorphism types of 
$n$-manifolds $M$ which admit Riemannian metrics satisfying (1.1). In 
addition, the convergence is in the weak $L^{2,p}$ topology, and the 
limit metric $g_{\infty}$ is $L^{2,p},$ for any $p <  \infty .$

 While conceptually important, this result is of somewhat limited applicability 
in itself, since the bound on the full curvature tensor $R$ is very 
strong and only realizable in special situations. A direct 
generalization of this result to the more natural situation where one 
imposes bounds on the Ricci curvature, i.e. 
\begin{equation} \label{e1.2}
|Ric| \leq  \lambda , \ vol \geq  v_{0}, \ diam \leq  D, 
\end{equation}
is false. For example, there are sequences of metrics on the double of 
the tangent bundle of $S^{2},$ i.e. $TS^{2}\cup_{\partial}TS^{2},$ 
satisfying (1.2), which converge in the Gromov-Hausdorff topology to a 
metric on the double of the cone on ${\Bbb R}{\Bbb P}^{3}, C({\Bbb 
R}{\Bbb P}^{3})\cup_{\partial}C({\Bbb R}{\Bbb P}^{3}),$ cf. [An2]. The 
limit is not a smooth manifold, but is an orbifold, with two orbifold 
singular points, (the vertices of the cones).

 On the other hand, in dimension 4, it was shown in [An3] that such 
degeneration to orbifolds (with point singularities) is the only 
possible degeneration under the bounds (1.2); any sequence $g_{i}$ of 
metrics on a fixed 4-manifold $M^{4}$ satisfying the bounds (1.2) has a 
subsequence converging in the Gromov-Hausdorff topology to an orbifold, 
with a bounded number of singular points. The convergence is 
$C^{1,\alpha}$ away from singular points. In fact, if one enlarges the 
class of manifolds to include orbifolds, then the space of orbifold 
metrics is compact in the $C^{1,\alpha}$ and weak $L^{2,p}$ topologies 
in a natural sense. These orbifold compactness results were first 
proved for Einstein metrics in [An1] and [BKN], and later generalized to 
the setting (1.2). Analogous results hold in all dimensions, and with 
varying manifolds, if one adds a fixed bound on the $L^{n/2}$ norm of 
the curvature tensor, cf. [An3].

\medskip

 In this paper, we consider orbifold compactness results which replace 
the apriori bound on the Ricci curvature in (1.2) with other curvature 
bounds. These results are then applied to prove orbifold compactness 
theorems for metrics which are critical points of other natural 
functionals on the space of metrics besides Einstein metrics. 

 To describe the results, a precise definition of the orbifolds that 
arise in this context is needed.

\medskip
\noindent
{\bf Definition.}
 An orbifold $V$ is a topological space which is a smooth manifold off 
a finite set of singular points $\{q_{k}\}.$ A neighborhood of each 
singular point is a finite union of cones on spherical space-forms 
$C(S^{n-1}/\Gamma)$, $\Gamma  \subset SO(n)$, where the vertex of each 
cone is identified with the point $q_{k}$. A metric $g$ on $V$ is 
$C^{m,\alpha}$ if $g$ is $C^{m,\alpha}$ (locally) on $V \setminus 
\{q_{k}\}$, and, in a local uniformization $B^{n}\setminus \{0\}$ of 
each cone, the metric $g$ extends to a $C^{0}$ metric on the ball 
$B^{n}$.

\medskip

 On the one hand, this definition is more restrictive than the 
definition due to Thurston [Th], which allows for codimension 2 
cone-like singular sets. On the other hand, the definition allows for 
several spherical cones to be joined at a single vertex. An orbifold 
will be called irreducible if there is exactly one cone at each 
singular point.

\medskip

 We first describe the main result in general, and then discuss 
applications to critical metrics. Let $B_{x}(r)$ denote the geodesic 
ball of radius $r$ about $x$ and $S_{x}(r)$ the corresponding geodesic 
sphere in a complete Riemannian manifold $(M, g)$. For simplicity, all 
manifolds are assumed to be connected and oriented.
\begin{theorem} \label{t 1.1.}
  Let ${\cal M}  = {\cal M} (\Lambda , \nu_{0}, \delta_{0}, C_{0})$ be 
the space of closed $n$-dimensional Riemannian manifolds $(M, g)$, $n \geq 3,$ 
satisfying the bounds
\begin{equation} \label{e1.3}
vol_{g}M = 1, \int_{M}|R|^{n/2}dV_{g} \leq  \Lambda , 
\end{equation}
together with the following two conditions:

 (i). (Non-Collapse). There exists $\nu_{0} > $ 0 such that 
\begin{equation} \label{e1.4}
vol B_{x}(r) \geq  \nu_{0}r^{n}, 
\end{equation}
for all $x\in M$ and $r \leq $ 1.

 (ii) (Small curvature estimate). There exists $\delta_{0} > $ 0 and 
$C_{0} <  \infty $ such that if $(\int_{B_{x}(r)}|R|^{n/2})^{2/n} \leq  
\delta_{0},$ then
\begin{equation} \label{e1.5}
\sup_{B_{x}(r/2)}|R| \leq  (\frac{C_{0}}{vol 
B_{x}(r)}\int_{B_{x}(r)}|R|^{n/2}dV_{g})^{2/n}. 
\end{equation}

 Then for any fixed values of $\Lambda , \nu_{0}, \delta_{0}, C_{0}$ 
and $\alpha  <$ 1, the space ${\cal M} $ is $C^{1,\alpha}$ orbifold 
compact, in that any sequence of metrics $g_{i}$ on n-manifolds $M_{i}$ 
satisying the bounds (1.3)-(1.5) has a subsequence converging, modulo 
diffeomorphisms, to a limit orbifold (V, g), with $C^{1,\alpha}$ metric 
$g_{\infty}.$ The convergence $(M_{i}, g_{i}) \rightarrow $ (V, g) is 
in the Gromov-Hausdorff topology, and in the $C^{1,\alpha}$ topology, 
uniformly on compact subsets of $V_{0} = V\setminus\{q_{k}\},$ where 
$\{q_{k}\}$ is the finite set of singular points. In particular, there 
are smooth embeddings of $V_{0}$ into $M_{i},$ for $i$ sufficiently 
large. The limit metric $g_{\infty}$ is locally in $L^{2,p}$ and the 
sequence converges uniformly on compact sets in the weak $L^{2,p}$ 
topology on $V_{0},$ for any $p<\infty .$ 

 In addition, there is a bound on the number of orbifold singularities, 
the number of cones, and the order of local fundamental groups, 
depending only on $\Lambda , \nu_{0}, \delta_{0}$ and $C_{0}.$
\end{theorem}

 We make some comments on the hypotheses. The volume bound (1.3) is of 
course just a normalization of $(M, g)$. The bound (1.4) requires that 
the volume of balls is not too small compared with the volume of 
Euclidean balls of the same radius. Without such a hypothesis, the 
structure of the singular set of limits is much more complicated, and 
has not been analysed in detail. The local curvature estimate (1.5), as 
well as the global bound (1.3) are analogues of the curvature bound 
(1.1) in Gromov's compactness theorem. Both hypotheses can be realized 
for natural classes of metrics, especially on 4-manifolds. On the other 
hand, in analogy to (1.2), the bound on $\sup_{B(r/2)} |R|$ in (1.5) 
can be weakened to a bound on $\sup_{B(r/2)} |Ric|$, cf. Remark 2.11. 

  It is a consequence of the proof of Theorem 1.1 that there are only  
finitely many diffeomorphism types of $n$-manifolds in ${\cal M}(\Lambda, 
\nu_{0}, \delta_{0}, C_{0})$, cf. Remark 2.10. Also, Theorem 1.1 holds if 
${\cal M}$ is replaced by the corresponding space of orbifolds, cf. Remark 2.12. 

 There are numerous examples of metrics satisfying the bounds 
(1.3)-(1.5) which do in fact converge to such orbifold limits; see \S 
2.3 and \S 3 for further discussion. Although the regular part $V_{0}$ 
of the limit orbifold $V$ embeds in $M_{i}$, for $i$ sufficiently 
large, in general the topology of the complement $M \setminus V_{0}$ 
can be arbitrary, at least if $\Lambda$ in (1.3) is sufficiently large, cf. 
Remark 2.9. Note that when $n = dim M$ is odd, there are no orientable 
$(n-1)$-dimensional spherical space forms except $S^{n-1}$. Hence, all 
cones at singular points are balls $B^{n}$. Even if the limit $V$ is a 
smooth manifold $M_{\infty}$, and $M_{i}$ is diffeomorphic to 
$M_{\infty}$, the convergence of the metrics $g_{i}$ may not be smooth, 
i.e. $C^{1,\alpha}$, near the points $q_{k}$.

\medskip

 Suppose that one adds the bound
\begin{equation} \label{e1.6}
vol B_{x}(r) \leq  V_{0}r^{n}, \forall x\in M, r \leq  1, 
\end{equation}
for some $V_{0} <  \infty $ to the hypotheses of Theorem 1.1. Then the 
result follows rather easily from the proof of orbifold compactness 
under the bounds (1.2) in [An3]; the main point is that the bounds 
(1.2) imply the bound (1.6) via the Bishop-Gromov volume comparison 
theorem. The proof of Theorem 1.1 when (1.6) holds is given in \S 2.1.

 With this understood, the main issue in the proof of Theorem 1.1 is to 
prove that (1.6) follows from the bounds (1.3)-(1.5). To do this, we 
need the following minor variation of standard definitions of 
asymptotically locally Euclidean spaces.

\medskip
\noindent
{\bf Definition.}
 A complete Riemannian $n$-manifold $(N, g)$ is ALE (asymptotically 
locally Euclidean) if there is compact set $K \subset  N$ such that 
each component $E_{j}$, $j = 1, ..., \kappa$ of $N\setminus K$ is 
diffeomorphic to $({\Bbb R}^{n}\setminus B)/\Gamma_{j},$ where $B$ is 
the unit ball about $\{0\}\in {\Bbb R}^{n}$ and $\Gamma_{j}$ is a 
finite subgroup of $SO(n)$, acting isometrically on ${\Bbb R}^{n}$. 
Further, there is a chart $\Phi : {\widetilde E_{j}} \rightarrow  {\Bbb 
R}^{n} \setminus B$ in which the metric $\widetilde g = \Phi^{*}g$ 
satisfies
\begin{eqnarray} \label{e1.7}
\widetilde g_{ab} = \delta_{ab} + o(1), \\ 
|\partial\widetilde g_{ab}|_{C^{\alpha}} = O(r^{-1-\alpha}), \nonumber
\end{eqnarray}
for all $\alpha  < 1$.

 The following characterization of ALE spaces is then the main tool 
used to establish (1.6).
\begin{theorem} \label{t 1.2.}
  Let (N, g) be a complete, non-compact Riemannian n-manifold, $n \geq 
$ 3, satisfying
\begin{equation} \label{e1.8}
vol B_{x}(r) \geq  \nu_{0}r^{n}, 
\end{equation}
for some $\nu_{0} > $ 0 and all $x\in N, r <  \infty .$ Suppose also 
that
\begin{equation} \label{e1.9}
\int_{N}|R|^{n/2}dV \leq  \Lambda  <  \infty , 
\end{equation}
and that the small curvature estimate (1.5) holds, for some given 
$\delta_{0}, C_{0}.$

 Then (N, g) is ALE, and the number $\kappa $ of ends is bounded by 
$\nu_{0}, \Lambda , \delta_{0}$ and $C_{0}.$ Further, there is a 
constant $V_{0},$ depending only on $\nu_{0}, \Lambda , \delta_{0}$ and 
$C_{0},$ such that for all $x$ and r,
\begin{equation} \label{e1.10}
vol B_{x}(r) \leq  V_{0}r^{n}. 
\end{equation}
The asymptotic volume ratio is given by
\begin{equation} \label{e1.11}
\lim_{r\rightarrow\infty}\frac{vol B_{x}(r)}{r^{n}} = 
\sum_{j=1}^{\kappa}\frac{\omega_{n}}{|\Gamma_{j}|}, 
\end{equation}
where $\omega_{n}$ is the volume of the unit ball in ${\Bbb R}^{n}.$
\end{theorem}

 The proof of Theorem 1.2 is given in \S 2.2. The fact that Theorem 1.2 
implies (1.6) and thus the remainder of the proof of Theorem 1.1 is 
given in \S 2.3. 

\medskip

 Next we state some applications of Theorem 1.1 to special classes of 
Riemannian metrics, particularly on 4-manifolds. As perhaps the main 
example, consider the $L^{2}$ norm of the Weyl curvature as a 
functional on the space of metrics on a 4-manifold $M = M^{4}:$
\begin{equation} \label{e1.12}
{\cal W} (g) = \int_{M}|W_{g}|^{2}dV_{g}. 
\end{equation}
Critical points of this functional are Bach-flat metrics, i.e. metrics 
satisfying the Euler-Lagrange equation
$$\delta\delta W + \tfrac{1}{2}W(Ric) = 0,$$
where $W(\cdot  )$ denotes the action of $W$ on symmetric bilinear 
forms. Via the Bianchi identity, this is equivalent to 
\begin{equation} \label{e1.13}
\delta d(Ric - \frac{s}{6}g) - W(Ric) = 0. 
\end{equation}
The functional ${\cal W} $ is conformally invariant, and so one needs 
to choose a representative of the conformal class. For the purposes 
here, it is most useful to choose a Yamabe representative, i.e. a 
metric $\gamma $ of constant scalar curvature minimizing the Yamabe 
functional in conformal class $[\gamma ];$ (the fact that $\gamma $ may 
not be unique is of no consequence).

 Well-known examples of Bach-flat metrics on 4-manifolds are metrics 
conformal to Einstein metrics, conformally flat metrics, and 
half-conformally flat metrics, i.e self-dual or anti-selfdual metrics, 
for which $W^{-}$ or $W^{+} =$ 0, respectively. The analogue of Theorem 
1.1 in this context is the following:
\begin{theorem} \label{t 1.3.}
  Let $M$ be a fixed closed 4-manifold. The space ${\cal C}^{+}$ of 
Bach-flat unit volume Yamabe metrics on M, for which the scalar 
curvature $s$ of $g$ satisfies
\begin{equation} \label{e1.14}
s_{g} \geq  s_{0} > 0, 
\end{equation}
and the $L^{2}$ norm of the Weyl curvature satisfies
\begin{equation} \label{e1.15}
\int_{M}|W|^{2} \leq  \Lambda , 
\end{equation}
for some $s_{0} >$ 0 and $\Lambda  <  \infty ,$ is orbifold precompact 
in the $C^{\infty}$ topology.
\end{theorem}

 We note that since Bach-flat metrics are critical for ${\cal W} ,$ the 
functional ${\cal W} $ is constant on connected components of moduli 
spaces of Bach-flat metrics. This result is the exact analogue 
of a result proved for Einstein metrics on 4-manifolds in [An1], [BKN]. 
Further discussion regarding Theorem 1.3, and related results for 
critical points of other natural curvature functionals on the space of 
metrics are proved in \S 3. For a selection of recent work related to the 
issues of this paper, see for instance, [Ak], [CH], [CQY], [Lo], [PT], and 
[TV1,2].

\medskip

 Finally we should include a comment on the history of this work. In analogy 
to the results established for Einstein metrics, G. Tian suggested to 
the author in 1993 that Theorem 1.3 might be true, and proposed working 
jointly on this project. Rather quickly, we realized that the main 
obstacle to carrying this out was to establish some version of Theorem 
1.2. By early 1995, the author had developed the basic ideas leading to 
a proof of (a version of) Theorem 1.2 and we had agreed that the arguments  
were in place to begin and complete a manuscript on the work. Unfortunately, 
both Tian and the author became involved with other, more pressing, 
projects and the work just lingered without further attention.

 With the recent appearance of [TV1], the author decided it was time to 
complete his ideas and approach to the problem, resulting in this 
paper. The results of [TV1] and [TV2] are related, but not identical to 
the results proved here.

 While the overall strategy of the proof of Theorem 1.2, (and the other 
results above), is the same as that developed by the author in 1994-95, 
more recent results on the structure of the cut locus of Riemannian 
manifolds have been incorporated, allowing for a simplification of certain 
technical aspects of the proof. 

 I would like to thank Gang Tian for suggesting that Theorem 1.3 might 
be true, and for our early collaboration on the project. My thanks also 
to John Lott and Marc Herzlich for remarks and correspondence which 
helped to clarify the exposition of the paper. 

\section{General Orbifold Compactness.}
\setcounter{equation}{0}

 This section is concerned with the proofs of Theorems 1.1 and 1.2. In 
\S 2.1, we discuss the proof of Theorem 1.1 in case the bound (1.6) 
holds, i.e.
\begin{equation} \label{e2.1}
vol B_{x}(r) \leq  V_{0}r^{n}, r \leq  1, 
\end{equation}
for some (arbitrary) constant $V_{0} <  \infty .$ The proof in this 
case is essentially identical to that of [An1], to which we refer for 
some further details, if needed. Following this, we explain which parts 
of this proof generalize to the case where (2.1) is not assumed, and 
indicate why Theorem 1.1 should then follow from Theorem 1.2. Theorem 
1.2 is proved in \S 2.2, while the proof of Theorem 1.1 is then quite 
easily completed in \S 2.3. 

\medskip

{\bf \S 2.1.}
 In this subsection, we prove Theorem 1.1 in case (2.1) holds. Thus, 
let $(M_{i}, g_{i})$ be any sequence of closed Riemannian $n$-manifolds 
satisfying (1.3)-(1.5) and (2.1). The volume bounds (1.4) and (2.1) 
imply a uniform upper bound $N(r/2,r)$ on the number of disjoint 
geodesic balls of radius $r/2$ in any $r$-ball, depending only on 
$\nu_{0}, V_{0}.$ Any maximal collection of disjoint balls of radius 
$1/2$ in $(M_{i}, g_{i})$ gives a covering of $(M_{i}, g_{i})$ 
by balls of radius 1. Hence, the volume bounds imply an upper bound on 
the cardinality of such a collection of balls, and consequently, one 
has a uniform upper bound on the diameter:
\begin{equation} \label{e2.2}
diam (M_{i}, g_{i}) \leq  D = D(\nu_{0},V_{0}). 
\end{equation}
Since one also has a uniform bound on the number $N(\varepsilon ,R)$ of 
disjoint $\varepsilon$-balls in any $R$-ball of $(M_{i}, g_{i})$, Gromov's 
weak compactness theorem [G] implies that a subsequence of $(M_{i}, g_{i})$ 
converges in the Gromov-Hausdorff topology to a complete length space 
$(X, d_{\infty})$. Thus we need to examine the structure of $(X, g_{\infty})$, 
and improve the convergence to the limit. 

 For $\delta_{0}$ as in (1.5), and any given $r > 0$, let $\{x_{k}\}$ 
be a maximal $r/2$ separated set, (depending on $i$), in $(M_{i}, g_{i})$. 
Thus, the geodesic balls $B_{x_{k}}(\frac{r}{2})$ are disjoint while the balls 
$B_{x_{k}}(r)$ cover $M_{i}$. (The choice of $\{x_{k}\}$ is of course not 
unique, but this is of no concern here). Let
\begin{equation} \label{e2.3}
G_{i}^{r} = \cup\{B_{x_{k}}(r): \int_{B_{x_{k}}(2r)}|R|^{n/2}dV < 
\delta_{0}\}, 
\end{equation}
where $R = R_{g_{i}}$ is the curvature tensor, and similarly, let
\begin{equation} \label{e2.4}
B_{i}^{r} = \cup\{B_{x_{k}}(r): \int_{B_{x_{k}}(2r)}|R|^{n/2}dV \geq  
\delta_{0}\}. 
\end{equation}
All quantities here are with respect to $(M_{i}, g_{i})$. Hence for each $i$, 
$M_{i} = G_{i}^{r}\cup B_{i}^{r}$. Observe that there is uniform bound 
on the number $Q_{i}^{r}$ of $r$-balls in $B_{i}^{r}:$
\begin{equation} \label{e2.5}
Q_{i}^{r} \leq  \frac{\Lambda}{\delta_{0}}, 
\end{equation}
independent of $i, r$. Further, exactly as in the proof of (2.2), by 
(1.4) and (2.1) there is a uniform bound on the intrinsic diameter of 
the geodesic annuli $A_{x}(\frac{1}{2}r, 2r)$, 
\begin{equation} \label{e2.6}
diam A_{x}(\frac{1}{2}r, 2r) \leq  Dr, 
\end{equation}
at arbitrary base points $x\in (M_{i}, g_{i})$, with $D = D(\nu , V_{0})$. 
For the same reasons, there is a uniform upper bound on the number of 
components of $A_{x}(\frac{1}{2}r,2r)$. By (2.1) and (2.5), the total volume 
of $B_{i}^{r}$ satisfies $volB_{i}^{r} \leq  V_{0}Q(2r)^{n}$, and so is 
small when $r$ is small. Hence, most of the volume of $(M_{i}, g_{i})$ 
is contained in $G_{i}^{r}.$

 The small curvature estimate (1.5) implies a uniform curvature bound 
on $G_{i}^{r}:$
\begin{equation} \label{e2.7}
|R| \leq  Cr^{-2} \ \ {\rm on} \ \  G_{i}^{r}. 
\end{equation}
In fact (2.7) holds on the $r/2$ thickening of $G_{i}^{r}.$ It then 
follows from a local version of the smooth (i.e. $C^{1,\alpha}\cap 
L^{2,p})$ Gromov compactness theorem, cf. [An1] or [An3] for example, 
that for any given $r > $ 0, a subsequence of $G_{i}^{r}$ converges in 
the $C^{1,\alpha}$ and weak $L^{2,p}$ topologies to a limit manifold 
$G_{\infty}^{r}$ with limit $C^{1,\alpha}\cap L^{2,p}$ metric 
$g_{\infty}^{r}.$ In particular, $G_{\infty}^{r}$ and $G_{i}^{r}$ are 
diffeomorphic, for $i$ large, and there exist smooth embeddings 
$F_{i}^{r}: G_{\infty}^{r} \rightarrow  G_{i}^{r} \subset  M_{i}$ such 
that $(F_{i}^{r})^{*}(g_{i})$ converges in $C^{1,\alpha}$ to 
$g_{\infty}^{r}.$

 Now choose a sequence $r_{j} \rightarrow $ 0, with $r_{j+1} = 
\frac{1}{2}r_{j},$ and perform the above construction for each $j$. Let 
$G_{i}(r_{m}) = \{x\in (M_{i}, g_{i}): x\in G_{i}^{j}, {\rm for \ some} 
\ j \leq  m\}$, so that one has inclusions
$$G_{i}(r_{1}) \subset G_{i}(r_{2}) \subset  ... \subset  M_{i} $$
By the argument above, each $G_{i}(r_{m}) \subset  (M_{i}, g_{i}),$ for 
$m$ fixed, has a subsequence converging smoothly, i.e. in $C^{1,\alpha}$ 
and weak $L^{2,p},$ to a limit $G_{\infty}(r_{m}).$ Clearly 
$G_{\infty}(r_{m}) \subset  G_{\infty}(r_{m+1})$ and we set
\begin{equation} \label{e2.8}
G = \cup_{1}^{\infty}G(r_{m}), 
\end{equation}
with the induced metric $g_{\infty}.$ Thus, $(G_{\infty}, g_{\infty})$ 
is $C^{1,\alpha}$ and $L^{2,p}$ smooth and for any $m$, there are smooth 
embeddings $F_{i}^{m}: G(r_{m}) \rightarrow  M_{i},$ for $i$ 
sufficiently large, such that $(F_{i}^{m})^{*}(g_{i})$ converges in 
$C^{1,\alpha}$ to the metric $g_{\infty}.$

 Let $\bar G$ be the metric completion of $G$ w.r.t. $g_{\infty}.$ Then 
there is a finite set of points $q_{k}$, $k = 1, ..., Q$, such that
\begin{equation} \label{e2.9}
\bar G = G \cup  \{q_{k}\}. 
\end{equation}
This follows since there is a uniform upper bound (2.5) on the 
cardinality of $B_{i}^{r},$ for all $r$ small, and all $i$, independent 
of $r, i$. It is then easy to see that a subsequence of $(M_{i}, g_{i})$ 
converges to the length space $(\bar G, g_{\infty})$ in the 
Gromov-Hausdorff topology, so that $X = \bar G$. Moreover, by the 
construction of $G$, for $i$ sufficiently large, and for any compact domain 
$K \subset G$, there are embeddings $F_{i}: K \rightarrow  M_{i}$ such that 
$F_{i}^{*}g_{i}$ converges to $g_{\infty}$ in the $C^{1,\alpha}$ and weak 
$L^{2,p}$ topologies on $K$. Thus, $g_{\infty}$ is locally in $C^{1,\alpha}$ 
and $L^{2,p}$.

 It remains to prove that $\bar G$ is an orbifold, with orbifold 
singular points $\{q_{k}\},$ i.e. $\bar G = V$. This follows by an 
analysis of the tangent cone of the limit metric $g_{\infty}$ near each 
$q_{k},$ i.e. by a blow-up analysis. The curvature of $G$ is locally 
bounded in $L^{p},$ for any $p <  \infty .$ Further, by lower 
semi-continuity of the norm under weak convergence, the $L^{n/2}$ norm 
of the curvature on $(G, g_{\infty})$ is globally bounded:
\begin{equation} \label{e2.10}
\int_{G}|R|^{n/2}dV \leq  \Lambda . 
\end{equation}
In particular, for any $q = q_{k}\in \bar G,$
\begin{equation} \label{e2.11}
\int_{A_{q}(\frac{1}{2}r,2r)}|R|^{n/2}dV \rightarrow  0, \ \ {\rm as} \ \  
r \rightarrow  0. 
\end{equation}
This, together with the estimate (1.5) implies that on the limit 
$(G, g_{\infty})$,
\begin{equation} \label{e2.12}
(\frac{1}{r^{n}}\int_{A_{q}(\frac{1}{2}r,2r)}|R|^{p}dV)^{1/p} << 
r^{-2}. 
\end{equation}
Observe that (2.10)-(2.12) are scale-invariant.

 Let $s_{j} = 2^{-j},$ for $j$ large, and rescale the metric 
$g_{\infty}$ on $G$ near $q$ by $s_{j}^{-2},$ i.e. consider the metrics 
$\bar g_{j} = s_{j}^{-2}\cdot  g_{\infty}.$ The volume bounds (1.4) and 
(2.1), together with the curvature bound (2.12) imply, via the local 
$C^{1,\alpha}\cap L^{2,p}$ Gromov compactness theorem, that a 
subsequence of $(G\setminus q, \bar g_{j}, q)$ converges, modulo 
diffeomorphisms, and in the $C^{1,\alpha}$ and weak $L^{2,p}$ 
topologies, to a flat limit $(T^{\infty}, \bar g^{\infty}).$ Using the 
diameter and volume bounds (2.6) and (2.1), it is straightforward to 
show that $T^{\infty}$ has a bounded number of components and the 
metric completion $\bar T^{\infty}$ of each component of $T^{\infty}$ 
has a single isolated singularity \{0\}. Thus, $\bar T^{\infty}$ is a 
finite collection of complete flat manifolds joined at a single 
isolated singularity $\{0\}$. From this, it follows also easily that 
$\bar T^{\infty}$ is isometric to a union of flat cones 
$C(S^{n-1}/\Gamma_{j})$. By the smooth convergence, this (unique) 
structure on the limit is equivalent to the structure of $(G, g_{\infty})$ 
on small scales near the singular point $q$, and shows that 
$q$ is an orbifold singularity; see [An1] for further details.

\medskip

 The remaining parts of Theorem 1.1 are now also easily established. 
Thus, the uniform bound on $Q_{i}^{r}$ in (2.5) gives a uniform bound 
on the number $Q$ of singular points in $(V, g_{\infty})$. The lower 
bound (1.4) on the local volumes implies a lower bound on the local 
volumes of each spherical cone $C(S^{n-1}/\Gamma_{j}),$ which gives a 
uniform upper bound on the order of each local fundamental group 
$\Gamma_{j}.$ The upper volume bound (2.1) then gives a uniform bound 
on the number of cones joined at any singular point $q \in  V$.

 As shown above, the limit metric $g_{\infty}$ is locally 
$C^{1,\alpha}$ and $L^{2,p}$ on $V_{0} = V\setminus\{q_{k}\}.$ A punctured 
neighborhood of any $q = q_{k}$ has a finite cover diffeomorphic to a 
finite collection of punctured balls $B^{n}\setminus\{0\}.$ The tangent 
cone analysis above shows that the metric becomes Euclidean to 
$0^{\rm th}$ order at \{0\}, so that the metric $g_{\infty}$ has a $C^{0}$ 
extension over \{0\} on each ball. (It is not to be expected in general that 
$g_{\infty}$ has a smooth extension across \{0\}).

 This completes the proof of Theorem 1.1 in case (2.1) holds.
{\endproof}

\begin{remark} \label{r 2.1. (i).}
  {\rm If one has a lower bound on the Ricci curvature of $(M_{i}, 
g_{i})$, $Ric_{g_{i}} \geq  \lambda g_{i}$, for some $\lambda  >  
-\infty$, then the limit orbifold $(V, g_{\infty})$ is irreducible, so 
that $V\setminus q_{k}$ is locally connected for all $k$. This is 
proved in [An1] by means of the Cheeger-Gromoll splitting theorem, 
[CGl]. Without such a lower bound on $Ric$, the number of cones may be 
greater than one, although it is necessarily finite.

{\bf (ii).}
 If $x_{i}\in (M_{i}, g_{i})$ is any (generic) sequence of points such 
that $x_{i} \rightarrow  q$ in the Gromov-Hausdorff topology, then the 
curvature of $(M_{i}, g_{i})$ blows up, i.e. diverges to infinity, as 
$i \rightarrow  \infty .$ If one rescales the metrics $g_{i}$ so that 
the curvature remains bounded near $x_{i},$ then a subsequence 
converges to a complete ALE space $(N, g^{\infty}),$ i.e. a complete 
manifold, or more generally, a complete orbifold, with a finite number 
of ends, each ALE. Blowing this limit $(N, g^{\infty})$ down, i.e. 
rescaling $g^{\infty}$ by factors converging to 0, gives a union of 
spherical cones joined at a single vertex \{0\}.

 Now the curvature of $(M_{i}, g_{i})$ may diverge to infinity at a 
number of different scales near any singular point $q$, giving rise to a 
collection of such ALE spaces associated with each scale. This gives 
rise to a so-called ``bubble-tree'' of ALE spaces and scales. The 
structure of the limit orbifold $(V, g_{\infty})$ near any singular 
point $q$ is recaptured by the structure at infinity of the complete 
ALE orbifold corresponding to the {\it  smallest}  rate at which the 
curvature of $(M_{i}, g_{i})$ diverges to infinity near $q$; this 
corresponds to the largest distance scale, since the curvature scale 
corresponds to the inverse square of the distance scale. For full 
details on the structure of such bubble trees and their relation with 
the orbifold limit, see [B] and [AC]; see also Remark 2.10 for further 
analysis. }
\end{remark}

 Consider now the proof of Theorem 1.1 without the assumption (2.1). In 
this circumstance, there is no uniform bound $N(r/2,r)$ on the number 
of $r/2$ balls in an $r$-ball. However, one can define the sets 
$G_{i}^{r}$ and $B_{i}^{r}$ in (2.3)-(2.4) in exactly the same way, and 
again the cardinality of $\{B_{i}^{r}\}$ is uniformly bounded, i.e. 
(2.5) holds. Hence, as $r \rightarrow $ 0, $\{B_{i}^{r}\}$ converges to 
finite number of points, (in a subsequence). On $G_{i}^{r},$ one has 
uniform curvature bounds (2.7), depending only on $r$. Standard 
comparison geometry then gives a uniform upper bound on the local 
volume ratios of geodesic balls contained in $G_{i}^{r}$. Hence, one 
can pass to a limit $G$, (in a subsequence), as in (2.8) just as before. 

 As in (2.9), let $\bar G$ be the metric closure of $G$. We note however 
that without the bound (2.1), it is apriori possible that $\bar G = G 
=$ \{pt\}; it could happen that all of the volume and curvature may 
concentrate in very small regions about a given point, see for example [Ak].

 We will prove that the bound (2.1) does in fact hold on $(M_{i}, 
g_{i})$ by using a blow-up rescaling argument, as discussed in Remark 
2.1(ii) with regard to the structure of the singularities of $V$. To do 
this however requires the characterization of ALE spaces given by 
Theorem 1.2. 

\medskip

{\bf \S 2.2.}
 This subsection is concerned with the proof of Theorem 1.2. Just as 
with Theorem 1.1, Theorem 1.2 is comparatively easy to prove if one 
assumes an upper bound on the volume growth as in (2.1), i.e.
\begin{equation} \label{e2.13}
vol B_{x_{0}}(r) \leq  {\cal V} r^{n}, 
\end{equation}
for some fixed base point $x_{0}$ and ${\cal V}  <  \infty ,$ cf. 
[An1,Thm.3.5]. On the other hand, the proof to follow is essentially a 
generalization of the proof in case (2.13) holds.

 We first prove the following special case of Theorem 1.2. This serves 
a basis for the more general result.
\begin{theorem} \label{t 2.2.}
  Let (N, g) be a complete, non-compact smooth Riemannian n-manifold, 
$n \geq $ 3, satisfying
\begin{equation} \label{e2.14}
vol B_{x}(r) \geq  \nu_{o}r^{n}, 
\end{equation}
for some $\nu_{o} > 0$ and all $x\in N$. Suppose in addition that 
there is a point $x_{0}$ and a non-increasing function $\varepsilon 
(r):[0,\infty ) \rightarrow  [0,\infty)$, with $\varepsilon (0) = 1$ 
and $\lim_{r\rightarrow\infty}\varepsilon (r) = 0$, such that
\begin{equation} \label{e2.15}
|R|(r) \leq  \frac{\varepsilon (r)}{1+r^{2}}, 
\end{equation}
where $|R|(r) = \sup_{x\in S_{x_{0}}(r)}|R|(x).$

 Then (N, g) is ALE. Further, the number $\kappa $ of ends of $N$ and 
the maximal volume ratio $V_{0} = \sup_{r}\frac{vol B_{x}(r)}{r^{n}}$ 
are bounded by the function $\varepsilon (r).$ The asymptotic volume 
ratio is given by
\begin{equation} \label{e2.16}
\lim_{r\rightarrow\infty}\frac{vol B_{x}(r)}{r^{n}} = 
\sum_{k=1}^{\kappa}\frac{\omega_{n}}{|\Gamma_{k}|}, 
\end{equation}
where $\omega_{n}$ is the volume of the unit ball in ${\Bbb R}^{n}.$
\end{theorem}

\begin{remark} \label{r 2.3.}
  {\rm Theorem 2.2 does not hold in dimension 2. Consider for example 
the following warped product metric on ${\Bbb R}^{2}:$
$$g = dr^{2} + f^{2}(r)d\theta^{2}, $$
where, for large $r$, $f(r) = r\ln r$, (for example). A simple calculation 
shows that the curvature of $g$ satisfies (2.15), while the volume 
growth satisfies (2.14). However, this metric is not ALE. In terms of 
the proof to follow below, the tangent cone at infinity is of the form 
of the universal cover $\widetilde{({\Bbb R}^{2}\setminus 0)}$ of 
$({\Bbb R}^{2}\setminus 0)$, which is not a cone on a spherical space-form; 
(it could be viewed as a cone with infinite cone angle). This situation is 
possible because $({\Bbb R}^{2} \setminus \{0\})$ is not simply connected. 
More importantly in the context of the proof below, the circle $S^{1}$ 
has no intrinsic curvature, whereas spheres $S^{n-1}$ have non-trivial 
intrinsic curvature when $n \geq $ 3. Note that in higher dimensions, 
one can form flat manifolds by taking products, 
$\widetilde{({\Bbb R}^{2}\setminus 0)} \times {\Bbb R}^{n-2}$ for example. }
\end{remark}

{\bf Proof of Theorem 2.2:}
 For clarity, we divide the proof into several steps.

\medskip

{\bf Step 1. (Comparison Geometry Estimates).}

 We begin with some standard estimates from comparison geometry. Fix 
the base point $x_{0},$ and let $r(x) = dist(x, x_{0}).$ Let $S(r)$ 
denote the geodesic sphere of radius $r$ about $x_{0}.$ By using 
comparison with a rotationally symmetric metric whose curvature 
satisfies the bound (2.15), one has the following well-known Hessian 
and Laplacian comparison estimates for $r$:
\begin{eqnarray} \label{e2.17}
D^{2}r \leq  \frac{1+\mu (r) }{r}\hat g, \\ 
\Delta r \leq  \frac{n-1+\mu (r)}{r} , 
\end{eqnarray}
where $\mu (r)$ depends only on $\varepsilon (r)$ and $\mu (r) 
\rightarrow 0$ as $r \rightarrow  \infty$. Here $D^{2}r$ is the 
Hessian and $\Delta r = tr D^{2}r$ is the Laplacian. The distance 
function $r$ is globally Lipschitz, with Lipschitz constant 1. However, 
it is not smooth on the cutlocus of $x_{0}.$ Off the cutlocus 
(2.17)-(2.18) hold in the usual sense, since the metric $g$ and 
function $r$ are smooth there. At such smooth points, $D^{2}r$ is the 
$2^{\rm nd}$ fundamental form $A$ of the geodesic sphere $S(r)$, while 
$\Delta r$ is the mean curvature $H$ of $S(r)$. However, the estimates 
(2.17)-(2.18) hold globally when $D^{2}r$ and $\Delta r$ are understood 
in a weak sense, either in the sense of support functions or as distributions; 
see for instance [Ka] for full details. Here $\hat g$ in (2.17) is 
given by $\hat g(X,X) = |X|^{2} - \langle X, \sigma' \rangle^{2}$, where 
$\sigma$ is any geodesic from $x_{0}$ to $\sigma (r)$; see also [BL] for another 
definition of $D^{2}r$ on most of the cut locus. (The exact form of 
(2.17) on the cut locus will not be used however). Well-known analogues 
of (2.18) appear in the proof of the Cheeger-Gromoll splitting theorem 
[CGl] and Gromov's extension of the Bishop volume comparsion theorem 
past the cut locus, [G].

 We also point out the obvious but important fact that the products 
$rD^{2}r$ and $r\Delta r$ are scale-invariant.

 For $r > 0$, let $v(r) = volS(r) > 0$, where $vol$ denotes the $(n-1)$-dimensional 
Hausdorff measure ${\cal H}^{n-1}$ induced by the metric $g$. 
\begin{lemma} \label{l 2.4.}
  The function v(r) is well-defined, for all $r > $ 0. If $r_{i}$ is an 
increasing sequence with $r_{i} \rightarrow  r <  \infty ,$ then
\begin{equation} \label{e2.19}
v(r) \leq  \liminf_{i\rightarrow\infty}v(r_{i}). 
\end{equation}
Moreover, let $\frac{d}{dr}^{+}v(r) = \limsup_{\varepsilon 
>0,\varepsilon\rightarrow 0}\frac{v(r+\varepsilon 
)-v(r)}{\varepsilon}.$ Then
\begin{equation} \label{e2.20}
\frac{d}{dr}^{+}v(r) \leq  \frac{n-1+\mu (r)}{r}v(r). 
\end{equation}
\end{lemma}

{\bf Proof:}
 Let ${\cal C}$ denote the cutlocus of $x_{0}$ in $(M, g)$. By [BL], 
off a set $X$ with ${\cal H}^{n-1}(X) = 0$, ${\cal C}$ is locally a 
smooth hypersurface, i.e. every point $x\in{\cal C} \setminus X$ has a 
neighborhood $N_{x}$ such that $N_{x}\cap{\cal C}  = H_{x}$, where 
$H_{x}$ is a smooth hypersurface in $M$. Points $y\in H_{x}$ are 
characterized by the fact that there exist exactly two minimal 
geodesics from $x_{0}$ to $y$, neither of which has a conjugate point 
at $y$. Since ${\cal H}^{n-1}(X) = 0$, $X$ will be ignored in the 
following. Let $H$ be the hypersurface given by the union of the local 
$(n-1)$-manifolds $H_{x}$ and let $S^{reg}(r)$ denote the regular set 
of $S(r)$, so that $S(r) = S^{reg}(r)\cup (H\cap S(r))$ up to a set of 
${\cal H}^{n-1}$-measure 0. By [IT, Thm. B], ${\cal H}^{n-1}({\cal C})$ 
is finite on compact subsets of $M$, and hence the same holds for the 
(relatively) open subset $H \subset  {\cal C}$. Since $H\cap S(r)$ is a 
closed and bounded subset of $H$, ${\cal H}^{n-1}(H\cap S(r))$ is 
finite, for all $r$. Hence $v(r)$ is defined for all $r$. Moreover, for 
all $r$ except a closed set of Hausdorff dimension $0$ in ${\Bbb R}$, 
$H\cap S(r)$ is of ${\cal H}^{n-1}$-measure 0; for such generic $r$, 
\begin{equation} \label{e2.21}
v(r) = volS^{reg}(r).
\end{equation}
An $r$ satisfying (2.21) will be called regular.

  The property (2.19) corresponds to a well-known property of the cutlocus. 
Since any $x\in S(r)$ is the endpoint of a collection of minimal geodesics 
$\sigma_{x}$ from $x_{0},$ and since such geodesics are unique up to 
their endpoints at any distance $r-\varepsilon$, $\varepsilon > 0$, it 
follows that $S(r)$ is the Hausdorff limit of $S^{reg}(r-\varepsilon )$ 
as $\varepsilon  \rightarrow $ 0. Any domain $U \subset  S^{reg}(r)$ is 
the smooth limit of (essentially unique) domains $U_{\varepsilon} 
\subset  S^{reg}(r-\varepsilon ),$ while for any given component $C$ of 
$H\cap S(r),$ there are exactly two disjoint domains 
$C^{\pm}(r-\varepsilon ) \subset  S^{reg}(r-\varepsilon )$ converging 
smoothly to $C$ as $\varepsilon  \rightarrow $ 0. In other words, the 
limit $C$ is doubly covered. These properties imply (2.19).

 To prove (2.20), suppose first that $r$ is non-regular, so that $H\cap 
S(r) \neq  \emptyset$, and hence $v(r) > vol S^{reg}(r)$. Then 
${\cal H}^{n-1}(H) < \infty$ implies that, for $\varepsilon$ 
sufficiently small, $v(r+\varepsilon )-v(r) \leq 0$. Thus, (2.20) 
holds in this case.

 Assuming then that $r$ is regular, we first observe that
\begin{equation} \label{e2.22}
\int_{A(r,r+\varepsilon )}\Delta r \geq  v(r+\varepsilon ) - v(r). 
\end{equation}
To see this, one applies the divergence theorem to a smooth domain 
$U_{\delta} \subset  A(r,r+\varepsilon)$ approximating 
$A(r,r+\varepsilon)$ with $U_{\delta}\cap{\cal C}  = \emptyset$. The 
boundary terms approximate integrals over the inner and outer boundaries 
$S(r+\varepsilon)$ and $S(r)$ of $A(r, r + \varepsilon)$ and the cut locus 
in the interior. The same argument as above establishing (2.19) shows that 
the outer boundary integral dominates $v(r+\varepsilon)$, (because of 
the double covering of $H\cap S(r)$), and that the boundary term in the 
interior approximating the cut locus is non-negative. (This argument is 
the same as that used in the proof of the splitting theorem for Ricci 
curvature, [CGl]). On the other hand, since $r$ is regular, (2.21) holds 
and the inner boundary integral converges to $v(r)$. 

 Using (2.22) together with (2.18) gives 
$$\frac{n-1+\mu (r)}{r}\sup_{[r,r+\varepsilon ]}v(s) \geq  
\frac{1}{\varepsilon}(v(r+\varepsilon )-v(r)). $$
Letting $\varepsilon  \rightarrow $ 0 and using the fact that $r$ is 
regular again gives (2.20).
{\endproof}

 By integrating (2.20) w.r.t. $r$, one obtains
\begin{equation} \label{e2.23}
\frac{v(r)}{r^{n-1}} \leq  Cr^{\delta (r)}, \ {\rm and} \ 
\frac{v(s)}{s^{n-1}} \leq 2^{\delta(r)}\frac{v(r)}{r^{n-1}}, s \in [r,2r],
\end{equation}
where $\delta (r) \rightarrow 0$ as $r \rightarrow  \infty$. In fact, 
rewriting (2.20) in the form $\frac{d}{dr}^{+}(v/r^{n-1}) \leq  
(v/r^{n-1})(\varepsilon (r)/r),$ it then follows that
\begin{equation} \label{e2.24}
\limsup_{r\rightarrow\infty}\frac{d}{dr}^{+}\frac{v(r)}{r^{n-1}} \leq  
0. 
\end{equation}

 These global statements on the volume growth are thus simple 
consequences of the infinitesimal estimate (2.18). We will also need 
corresponding local estimates on the volume growth. Thus, for any given 
$r_{0} > $ 1, let $D(r_{0})$ be any domain in $S(r_{0})$ and let $D(r) 
\subset S(r)$ be defined as
\begin{equation} \label{e2.25}
D(r) = \{x\in S(r): dist(x, D(r)\} \leq  r-r_{0}\}, 
\end{equation}
for $r \geq  r_{0}.$ If $r \leq  r_{0}$, $D(r)$ is defined in the same 
way, with $|r-r_{0}|$ in place of $r-r_{0}.$ Then the same argument as 
above gives
\begin{equation} \label{e2.26}
\limsup_{r\rightarrow\infty}\frac{d}{dr}^{+}\frac{volD(r)}{r^{n-1}} 
\leq  0. 
\end{equation}

 The estimates (2.17)-(2.18) and hence (2.24)-(2.26) require only a 
lower bound on the curvature, i.e. follow from the bound
$$(K_{\min})(r) \geq  - \frac{\mu (r)}{r^{2}}, $$
where $(K_{\min})(r)$ is the minimal value of the sectional curvatures 
of $(N, g)$ on $S(r)$. From the upper bound on the curvature, i.e. 
$(K_{\max})(r) \leq  \frac{\mu (r)}{r^{2}},$ one obtains the following: 
suppose $x\in S(r)$ is not in the cutlocus of $x_{0},$ so that there is 
a unique geodesic $\sigma $ from $x_{0}$ to $x$. The geodesic $\sigma $ 
has length $r$. Suppose moreover that $\sigma (t)\in S(t)$ is not in the 
cutlocus of $x_{0}$ for any $t \leq  2r$, (or $t \leq  (1+\delta )r,$ 
for some $\delta  > 0$). Then at such $x$,
\begin{eqnarray} \label{e2.27}
A = D^{2}r \geq  \frac{1-\mu (r) }{r}\hat g, \\
H = \Delta r \geq  \frac{n-1-\mu (r)}{r}, \nonumber
\end{eqnarray}
(where $\mu (r)$ also depends on $\delta)$. These estimates will 
however not hold, in any generalized sense, at the cutlocus.
\begin{remark} \label{r 2.5.}
  {\rm Suppose the manifold $(N, g)$ has a pole, i.e. a point, say 
$x_{0},$ at which the exponential map is a diffeomorphism. Then the 
geodesic spheres $S(r)$ are connected, smooth manifolds, for all $r$. 
Consider the rescaled Riemannian metrics $g_{r} = r^{-2}g_{S(r)}$ on 
$S(r)$. From (2.17), (2.27) and the Gauss equations for the hypersurface 
$S(r) \subset N$, one has 
$$|K_{g_{r}} - 1| = \mu (r) \rightarrow  0, \ \ {\rm as} \ \  
r \rightarrow  \infty , $$
where $K_{g_{r}}$ is the intrinsic sectional curvature of 
$(S(r), g_{r}).$

 Thus, the geodesic spheres, when rescaled to unit radius, have 
curvature uniformly converging to 1. Hence, there is a uniform upper 
bound on their volume and diameter, by Myer's theorem. Translating this 
statement back to the unscaled metric g, it follows that 
$\limsup_{r\rightarrow\infty} volS(r)/r^{n-1} \leq  \sigma_{n-1},$ where 
$\sigma_{n-1}$ is the volume of the Euclidean $(n-1)$-sphere of radius 
1. Integrating this inequality, it follows that (2.13) holds. In 
particular, it follows that $(N, g)$ is ALE, in fact $(N, g)$ is AE, 
(asymptotically Euclidean). Note that this proof requires $n \geq 3$; 
for $n = 2$ the spheres are circles, and so have no intrinsic 
curvature, as in Remark 2.3. Note also this argument does not require a 
lower bound on the volume growth (2.14).

 The proof of Theorem 2.2 is a generalization of this kind of argument 
to the situation where there are no poles.

 This simple argument is due to Gromov, (unpublished), and gives rather 
simple proofs of results of Siu-Yau [SY] and Greene-Wu [GW] for instance. }
\end{remark}

{\bf Step 2. Tangent cones at infinity.}

 Let $r_{k}$ be any sequence diverging to infinity, and consider the 
geodesic annuli $A_{k} = A_{k}(m) = A_{k}(m^{-1}r_{k}, mr_{k})$ about 
$x_{0},$ for some arbitrary but fixed 1 $<  m <  \infty .$ We rescale 
the metric $g$ on $A_{k}$ by setting 
\begin{equation} \label{e2.28}
g_{k} = r_{k}^{-2}\cdot  g. 
\end{equation}
All distances of $g$ are contracted by the factor $r_{k}^{-1}$ w.r.t. 
$g_{k},$ so that $(A_{k}, g_{k})$ is a geodesic annulus of width 
$(m^{-1}, m)$. The curvature bound (2.15) implies that, w.r.t. $g_{k},$
\begin{equation} \label{e2.29}
|R_{g_{k}}|(x) \leq  \varepsilon (r_{k}), 
\end{equation}
for all $x\in A_{k},$ with $\varepsilon (r_{k}) \rightarrow $ 0 as $k 
\rightarrow  \infty .$ Further, the lower volume bound (2.14) implies 
that
\begin{equation} \label{e2.30}
vol B_{x_{k}}(s) \geq  \nu_{o}s^{n}, 
\end{equation}
for all balls $B_{x_{k}}(s) \subset  (A_{k}, g_{k}).$

 It follows from (the local $C^{1,\alpha}\cap L^{2,p}$ version of) the 
Gromov compactness theorem, cf. [An1], [An3], that a subsequence of the 
Riemannian manifolds $(A_{k}, g_{k}, x_{k}),$ at any base point 
sequence $x_{k}$ with $r(x_{k}) = r_{k},$ converges, (modulo 
diffeomorphisms), in the $C^{1,\alpha}$ and weak $L^{2,p}$ topologies, 
to a $C^{1,\alpha}\cap L^{2,p}$ limit Riemannian manifold $(A_{\infty}, 
g_{\infty}, x_{\infty}).$ By definition, $A_{\infty}$ is connected, 
containing the base point $x_{\infty} = \lim x_{k}$ and of course 
$A_{\infty} = A_{\infty}(m).$ On the limit $(A_{\infty}, g_{\infty}),$ 
the curvature is well-defined in $L^{p}.$ Since the $L^{p}$ norm of the 
curvature is lower semi-continuous, (2.29) implies that the limit 
$(A_{\infty}, g_{\infty})$ is flat. In particular, the limit is 
$C^{\infty}$ smooth. 

 Next, choose a sequence $m_{j} \rightarrow  \infty ,$ and perform the 
analysis above for each $m_{j}.$ It follows that a diagonal subsequence 
of the double sequence $(A_{k}(m_{j}), g_{k}, x_{k})$ converges to a 
maximal connected limit $(T_{\infty}, g_{\infty}, x_{\infty}).$ Such 
limits $T_{\infty}$ are called {\it tangent cones at infinity}, (from 
terminology originating in the theory of minimal varieties), although 
apriori they need not be cones. The arguments above prove that all 
tangent cones at infinity of $(N, g)$ are flat.

 Each tangent cone $T_{\infty}$ has a distinguished distance function 
$\rho  = r_{\infty}$ given as follows. For $r(x)$ as above, let $r_{k}(x) 
= r_{k}^{-1}\cdot  r(x).$ Then $r_{k}(x_{k}) =$ 1, and $r_{k}$ is a 
distance function on $A_{k}.$ The functions $r_{k}$ are Lipschitz, with 
Lipschitz constant 1. In the subsequence above, the sequence 
$\{r_{k}\}$ converges, in the Lipschitz topology, to the limit distance 
function $r_{\infty} = \rho $ on $T_{\infty}.$ By construction, $\rho $ 
takes on all values in (0, $\infty )$ on $T_{\infty}.$

 The level sets $S(t) = \{x\in T_{\infty}: \rho (x) = t\}$ of $\rho$ 
are called the geodesic spheres of $T_{\infty}.$ Thus, $S(t)$ is obtained 
by taking a limit of the rescaling of $S(tr_{k})$ by $r_{k}^{-1}.$ 
Observe that in general only part of $(S(tr_{k}), g_{k})$ will limit on 
$S(t) \subset  T_{\infty}$; other parts may be at infinite distance to 
$T_{\infty}$. Further, $S(t)$ is not necessarily connected, and the 
number of components may change with $t$. Observe also that the base 
point $x_{\infty}$ of $T_{\infty}$ lies in $S(1)$.

\medskip

{\bf Step 3. Geometry of Tangent Cones.}

 All tangent cones $T_{\infty}$ are flat manifolds, possibly with 
singularities at the metric boundary $\partial T_{\infty}$. 
Neighborhoods of the boundary $\partial T_{\infty}$ are given by 
the balls $B(\varepsilon) = \{x: \rho(x) < \varepsilon\}$.  From 
standard scaling properties, the estimates (2.17)-(2.18) translate to
\begin{eqnarray} \label{e2.31}
D^{2}\rho  \leq  \frac{1}{\rho}\hat g, \\
\Delta\rho  \leq  \frac{n-1}{\rho}, 
\end{eqnarray}
where again (2.31)-(2.32) are understood to hold in the sense of 
support functions or distributions at cut points of $\rho .$ 

 In particular, for any bounded domain $D \subset S(1) \subset  
T_{\infty}$ and associated $D(t) \subset S(t)$ as in (2.25), one has, 
as in (2.20),
\begin{equation} \label{e2.33}
\frac{d}{dt}^{+}\frac{vol D(t)}{t^{n-1}} \leq  0, 
\end{equation}
By Lemma 2.4, for a bounded domain $D$, $vol D(t) \subset S(1)$ is finite.

 Now return for the moment to the original manifold $(N, g)$. For any 
given $R$ large, decompose $S(R)$ into a collection of connected 
subdomains as follows. Let $C_{x}(s) = \{y\in S(R): 
dist_{A(R-\varepsilon, R+\varepsilon)}(y, x) <  s\}$, where the distance 
is taken in the annulus $A(R-\varepsilon, R+\varepsilon)$ about $x_{0}$; 
here $\varepsilon = \varepsilon(R) < 10^{-3}$ is taken sufficently small 
so that the number of components of $A(R-\varepsilon, R+\varepsilon)$ equals 
the number of components of $S(R)$. If $x$ and $x'$ are in different 
components of $S(R)$, (and hence in different components of 
$A(R-\varepsilon, R+\varepsilon)$), define as usual 
$dist_{A(R-\varepsilon,R+\varepsilon)}(x, x') = \infty$. 

 We may then cover $S(R)$ by a finite collection of connected closed sets 
$D_{i},$ with disjoint interiors, such that
\begin{equation} \label{e2.34}
C_{x_{i}}(2\pi R) \subset  D_{i} \subset  C_{x_{i}}(4\pi R). 
\end{equation}
for some collection of points $x_{i}\in S(R).$ To be precise, choose a 
maximal collection of disjoint balls $C_{x_{i}}(2\pi R)$ in $S(R)$ so 
that the balls $C_{x_{i}}(4\pi R)$ cover $S(R)$. Then choose the sets 
$D_{i}$, so that $C_{x_{i}}(2\pi R) \subset  D_{i}$ with $\partial 
D_{i}$ contained in $C_{x_{i}}(4\pi R)\setminus C_{x_{i}}(2\pi R)$ and 
${\cal H}^{n-1}(\partial D_{i}) = 0$. Let $N(R)$ denote the cardinality of 
the collection $\{D_{i}\}$. Of course it is possible apriori that $N(R)$ 
diverges to infinity as $R \rightarrow  \infty$.

 By construction, we then have
\begin{equation} \label{e2.35}
\frac{v(R)}{R^{n-1}} = \sum_{i=1}^{N(R)}\frac{vol(D_{i})}{R^{n-1}}. 
\end{equation}
Again, for any fixed $R$ large and $D_{i} = D_{i}(R)$ as above, define 
the domains $D_{i}(r) \subset S(r)$, as in (2.25), for $\frac{3}{4}R \leq  r 
\leq  \frac{5}{4}R$. Then 
$$S(r) = \cup_{i=1}^{N(R)}D_{i}(r),$$
and the structure of $S(r)$ discussed in the proof of Lemma 2.4 and the triangle 
inequality imply that the domains $D_{i}(r)$ satisfy 
${\cal H}^{n-1}(D_{i}(r) \cap D_{j}(r)) = 0$, for $i \neq j$. Observe that 
$D_{i}(r)$ may no longer be necessarily connected, for $r \neq R$, but this 
is of no importance. By construction, we then have
\begin{equation} \label{e2.36}
\frac{d}{dr}^{+}(\frac{v(r)}{r^{n-1}})_{r=R} = 
\frac{d}{dr}^{+}(\sum_{i=1}^{N(R)}\frac{vol(D_{i}(r))}{r^{n-1}})_{r=R} = 
\sum_{i=1}^{N(R)}\frac{d}{dr}^{+}(\frac{v_{i}(r)}{r^{n-1}})_{r=R}, 
\end{equation}
where $v_{i}(r) = vol D_{i}(r)$.

  Before proceeding further, we need a slight refinement of the considerations 
above. Let $S^{-}(R)$ be the union of the components of $S(R)$ which 
(partially) bound a domain in the annulus $A(R,2R)$ about $x_{0}$. Thus, 
a component $C(R) \subset S^{-}(R)$ if $C(R) \subset \partial K(R)$, where 
$K(R)$ is a domain in $A(R,2R)$ with $\partial K(R) \subset S(R)$. Let 
$S^{+}(R) = S(R) \setminus S^{-}(R)$, so that $S(R) = S^{+}(R) \stackrel{\cdot}{\cup} 
S^{-}(R)$. Note that $S^{+}(R)$ is always non-empty; it contains components which are 
the boundaries of ends of $N$ outside $S(R)$, while $S^{-}(R)$ may or may not 
be empty. Let $v^{+}(R) = vol S^{+}(R)$, $v^{-}(R) = vol S^{-}(R)$, so that 
$v(R) = v^{+}(R) + v^{-}(R)$, for all $R$. Observe that (2.35)-(2.36) hold with 
$v^{+}(r)$ or $v^{-}(r)$ in place of $v(r)$.

   Let $D^{-}(R) \subset A(R,2R)$ be the compact domain with $\partial D(R) 
= S^{-}(R)$. Then (2.23) gives
\begin{equation} \label{e2.37}
\frac{vol (S(s) \cap D^{-}(R))}{s^{n-1}} \leq (1+\delta)\frac{vol S^{-}(R)}{R^{n-1}},
\end{equation}
where $\delta \rightarrow 0$ as $R \rightarrow \infty$. Thus, the main issue is to 
control the volume growth of $S^{+}(R)$ as $R \rightarrow \infty$.

 Now suppose
\begin{equation} \label{e2.38}
\frac{d}{dr}^{+}\frac{v^{+}(r)}{r^{n-1}} \leq  0, 
\end{equation}
for all $r \geq  R_{0}$, for some $R_{0} <  \infty$. Then (2.37) and (2.38) 
clearly imply $v(r) \leq  Kr^{n-1}$, for some $K < \infty$ and so 
\begin{equation} \label{e2.39}
vol B_{x_{0}}(r) \leq  Kr^{n}.
\end{equation}
Hence, as noted prior to its statement, Theorem 2.2 follows in this 
case from [An1], (see also Step IV below for details on the proof). 
Thus, we may suppose that there exists some sequence $r_{k} \rightarrow 
 \infty$ such that
\begin{equation} \label{e2.40}
\frac{d}{dr}^{+}\frac{v^{+}(r)}{r^{n-1}}|_{r=r_{k}} >  0. 
\end{equation}
It follows from (2.36) and (2.40) that there exist $i = i(r_{k})$ and 
domains $D_{i} \subset S^{+}(r_{k})$, as in (2.34) such that for all $k$,
\begin{equation} \label{e2.41}
\frac{d}{dr}^{+}\frac{v_{i}(r)}{r^{n-1}}|_{r=r_{k}} >  0. 
\end{equation}
Let $x_{i}\in S^{+}(r_{k})$, $i = i(k)$, be the base points of $D_{i}$. 
By passing to a subsequence based at $\{x_{i}\}$ if necessary, 
one obtains convergence to a tangent cone $(T_{\infty}, g_{\infty}, 
x_{\infty})$ with associated distance function $\rho$ and level sets 
$S(t)$, and with $x_{\infty} = \lim x_{i} \in S(1)$. The convergence is in the 
$C^{1,\alpha}$ topology (on compact subsets), while the convergence of the 
distance functions is at least Lipschitz. Let $D_{\infty} \subset S(1)$ be the 
limit of the domains $D_{i}$. Then, by construction, $D_{\infty}$ is not the 
(partial) boundary of a domain in $A(1,2) \subset T_{\infty}$, so that $\rho$ assumes 
all values in $(\frac{1}{2}, \frac{3}{2})$ near $D_{\infty}$. In particular, $D_{\infty}$ 
locally disconnects $T_{\infty}$ and hence $vol D_{\infty} = {\cal H}^{n-1}(D_{\infty}) 
> 0$.

  Note that (2.41) implies that the domain $D_{i} \subset S(r_{k})$ is regular, 
cf. (2.21). By taking a small perturbation of $\{r_{k}\}$ if necessary, (relative 
to its own size), we may assume that $D_{\infty} \subset S(1) \subset T_{\infty}$ is 
also regular and (2.41) still holds. For if this were not the case, then the upper 
bound (2.20), cf. also (2.23), would imply again that (2.39) holds. Let 
$D_{\infty}(t) \subset S(t)$ be defined as in (2.25), so that $D_{\infty}(t)$ 
is the limit of the the corresponding domains in $(S(t), g_{k})$. Also, let 
$v_{\infty}(t) = vol D_{\infty}(t) = {\cal H}^{n-1}(D_{\infty}(t))$. 
The regularity of these domains and the $C^{1,\alpha}$ convergence to the tangent 
cone $T_{\infty}$ implies that, for $t$ close to $1$, 
$$
\frac{vol_{g}(D_{i}(tr_{k}))}{(tr_{k})^{n-1}} = \frac{vol_{g_{k}}D_{i}(t)}{t^{n-1}} 
\rightarrow \frac{vol D_{\infty}(t)}{t^{n-1}},$$
and 
\begin{equation} \label{e2.42}
r_{k}\frac{d}{dr}^{+}\frac{vol_{g}(D_{i}(r))}{r^{n-1}}|_{r=r_{k}} = 
\frac{d}{dt}^{+}\frac{vol_{g_{k}}(D_{i}(t))}{t^{n-1}}|_{t=1} \rightarrow  
\frac{d}{dt}^{+}\frac{v_{\infty}(t)}{t^{n-1}}|_{t=1}.
\end{equation}
Hence, from (2.41) one has
\begin{equation} \label{e2.43}
\frac{d}{dt}^{+}\frac{v_{\infty}(t)}{t^{n-1}}|_{t=1} \geq  0. 
\end{equation}

 On the other hand, from (2.33) one has
$$\frac{d}{dt}^{+}\frac{v_{\infty}(t)}{t^{n-1}}|_{t=1} \leq  0. $$
so that 
\begin{equation} \label{e2.44}
\frac{d}{dt}^{+}\frac{v_{\infty}(t)}{t^{n-1}}|_{t=1} = 0. 
\end{equation}
The pointwise estimate (2.32), (in the sense of support functions), together with 
(2.44) implies that $\Delta \rho \equiv (n-1)$ on $D_{\infty}.$ Moveover, since 
$D^{2}\rho \leq  \hat g$, and $\Delta \rho  = tr D^{2}\rho$, it also follows that 
$D^{2}\rho = \hat g$ on $D_{\infty}$. This implies first that the domain 
$D_{\infty} \subset S(1)$ is in fact smooth. As in Remark 2.5, the Gauss equations 
imply that $D_{\infty}$ has constant curvature $1$, and hence is a domain in a spherical 
space form $S^{n-1}/\Gamma$, (since $n \geq 3$). Moreover, by the 
construction in (2.34), $D_{\infty}$ contains an intrinsic geodesic 
ball of radius at least $2\pi$, and so it follows that 
\begin{equation} \label{e2.45}
D_{\infty} = S^{n-1}/\Gamma . 
\end{equation}
Since $\rho$ is also a distance function from $D_{\infty},$ it follows 
that for $t \leq 1$, $D_{\infty}(t)$ is the spherical space-form of 
radius $t$, and curvature $t^{-2}.$ Thus, the domain $B_{\infty} = 
\cup_{t\leq 1}D_{\infty}(t) \subset  T_{\infty}$ is a ball of radius 1 
in the standard cone $C(S^{n-1}/\Gamma ),$ with vertex $\{0\} \subset  
\{t = 0\}$.

 At this stage, it is possible that $S(1) \subset  T_{\infty}$ has 
other components besides $D_{\infty};$ of course each component 
satisfying (2.43) is the unit ball about the vertex in a standard cone 
$C(S^{n-1}/\Gamma )$. 

  Summarizing, the argument above proves that on any sequence of 
base points $x_{k},$ with $r_{k} = r(x_{k}) \rightarrow  \infty$, for 
which 
\begin{equation} \label{e2.46}
\liminf_{r\rightarrow\infty}\frac{d}{dr}^{+}\frac{vol(D_{k})}{r^{n-1}}|_{
r=r_{k}} \geq  0, 
\end{equation}
any limit $D_{\infty}$ is of the form (2.45). Moreover, the domain $D_{k}$ 
is a full component of $S(r_{k})$ for $k$ sufficiently large.

\medskip

{\bf Step 4. Completion of Proof.}

 The statement (2.46) implies the following: there is an $R < \infty$, 
sufficiently large, and $\delta > 0$ sufficiently small, both 
depending on $(N, g)$, such that for any $r \geq  R$, and any domain 
$D_{k} \subset S^{+}(r)$ as in (2.34), either
\begin{equation} \label{e2.47}
\frac{d}{dr}^{+}\frac{vol(D_{k})}{r^{n-1}} < -\delta , 
\end{equation}
or, 
\begin{equation} \label{e2.48}
\frac{vol(D_{k})}{r^{n-1}} \leq  (1+\delta 
)\frac{\sigma_{n-1}}{|\Gamma_{k}|}. 
\end{equation}
Thus, the volume ratio of the domain $D_{k}(r)$ is either strictly 
decreasing with increasing $r$, or the volume ratio is close to that of 
the corresponding spherical space-form. In the latter case, the geometry 
of $(D_{k}, g_{k})$ is $\delta$-close to that of a spherical 
space-form $S^{n-1}/\Gamma_{k}.$

 These facts imply that $volS(r)/r^{n-1}$ is uniformly bounded. In 
fact, for all $r \geq R$,
$$\frac{vol(S(r))}{r^{n-1}} \leq  2\frac{vol(S(R))}{R^{n-1}},  $$
and hence one has
\begin{equation} \label{e2.49}
vol B(r) \leq  V_{0}r^{n}, 
\end{equation}
for some $V_{0}$ depending on $(N, g)$. This shows that (2.13) holds, and 
hence from the methods of [An1], it follows that $(N, g)$ has a finite 
number of ends, each ALE.

\medskip

 For completeness, we give the proof of this. The bounds (2.14) and 
(2.49) imply that the geodesic sphere $S(1)$ in any tangent cone at 
infinity $T_{\infty}$ of $(N, g)$ has uniformly bounded volume, diameter 
and a uniformly bounded number of components. The proof of this is the 
same as that given in \S 2.1, since one has a uniform upper bound on the number of 
disjoint balls of radius $\frac{1}{2}$ in any 1-ball $(B(1), g_{r})$ 
about $x_{0},$ for $g_{r}$ as in (2.28). Given a tangent cone 
$(T_{\infty}, g_{\infty})$ obtained as a limit of rescalings by factors 
$r_{k}^{-2}$ as in (2.28), observe that the geodesic sphere $S(t) 
\subset  T_{\infty}$ is just the unit geodesic sphere $S(1)$ in the 
tangent cone $(T_{\infty}(t), g_{\infty}(t))$ obtained by rescaling by 
the factors $(tr_{k})^{-2}.$

 It follows that in any tangent cone $T_{\infty},$ and for any $t >$ 0, 
one has the uniform bound
\begin{equation} \label{e2.50}
diam S(t) \leq  D\cdot  t, 
\end{equation}
for some fixed $D$, compare with (2.6). Here $diam$ is the intrinsic 
diameter of each component of $S(t)$. For the same reasons, there is a 
uniform upper bound on the number of components of $S(t)$. As discussed 
near the end of \S 2.1, these estimates imply that the metric 
completion $\bar T_{\infty}$ of any $T_{\infty}$ is given by adding a 
bounded number of points $\{z_{j}\}$ to $T_{\infty},$ one for each 
component of $S(t)$, for $t$ small. Hence $\bar T_{\infty}$ is a complete 
flat manifold with a finite collection of isolated singularities. 

 By an analysis of the developing map for flat structures, one sees 
that $T_{\infty}$ is a standard flat cone $C(S^{n-1}/\Gamma)$ with 
a finite number of points removed, (including the vertex). Thus, only one 
point, the vertex say $\{0\} = \{z_{0}\}$, among the collection 
$\{z_{j}\}$ may be a singular point; all others have neighborhoods 
which are balls, and so are locally removable singularities. This result 
is proved in [An4,Thm.3.2], cf. also [AC], the idea being that isolated 
singularities of flat manifolds are necessarily orbifold singularities, 
and all flat orbifolds are good orbifolds in the sense of Thurston, i.e. 
are quotients of ${\Bbb R}^{n}$ by a group of isometries.

 Thus {\it all} tangent cones $T_{\infty}$ of $(N, g)$ are standard cones on 
spherical space-forms, possibly with a finite number of points removed.
Moreover, we now claim that the collection of removable 
singularities $\{z_{j}\}$, $j \geq 1$, is empty, so that all geodesic spheres 
$S(1)$ are standard spherical space-forms $S^{n-1}/\Gamma$. Namely, since the 
convergence of the rescalings $r^{-2}\cdot  g$ to $g_{\infty}$ is smooth, 
(i.e. $C^{1,\alpha}),$ it follows that for any $r$ sufficiently large, and any 
base point $x\in S(r)$, the component $(A(\frac{1}{2},1), g_{r}, x)$ of 
$(A(\frac{1}{2},1), g_{r})$ containing $x$ is close, in $C^{1,\alpha}$, 
to the standard flat metric on a standard annulus $A(\frac{1}{2}, 1)$ 
in $C(S^{n-1}/\Gamma )$. Npw recall again that if $(S(1), g_{r})$, 
(the rescaling of $(S(r), g)$ by $r^{-2}$), is close to $S(1) \subset T_{\infty}$ 
then for any given $\lambda$ large, $(S(\lambda), g_{r})$ is close to 
$S(\lambda) \subset T_{\infty}$, but is also close to $S(1) \subset T_{\infty}^{\lambda}$, 
where $T_{\infty}^{\lambda}$ is $T_{\infty}$ scaled down by the factor $\lambda^{-2}$. 
It follows that for all $r$ sufficiently large, $r \geq  R_{0}$, each component of 
$A(\frac{1}{2}r, r) \subset  (N, g)$ stabilizes, i.e. is independent of $r$; no 
bifurcations or changes in the topology are possible, since the geometry of each 
annular region remains very close to that of a standard spherical cone. In 
particular, for $r$ sufficiently large, the components of $S(r)$ can not bifurcate 
or merge and hence each geodesic sphere $S(1)$ in any tangent cone at infinity 
$T_{\infty}$ is connected.

 Exactly the same reasoning implies that $(N, g)$ has only a finite 
number of ends, and each end is asymptotic to a standard cone 
$C(S^{n-1}/\Gamma ).$ Scaling properties and the $C^{1,\alpha}$ 
convergence to the tangent cone at infinity imply the decay estimates 
(1.7). The formula (2.16) follows then immediately.

 Finally, we claim that the number $\kappa $ of ends of $(N, g)$, and the 
volume constant $V_{0}$ depend only on the function $\varepsilon (r).$ 
Since $\varepsilon (r)$ is fixed, for $r$ sufficiently large, 
$\varepsilon (r)$ is sufficiently small. Hence, for a given $\mu $ 
small, the rescaled curvature $r^{2}|R|(r)$ satisfies $r^{2}|R|(r) \leq 
 \mu ,$ for all $r \geq  R_{0} = R_{0}(\mu ).$ As above, this implies 
that the topology of $A(R_{0}, \infty )$ is that of a finite collection 
of truncated cones on $S^{n-1}/\Gamma .$ On the other hand, for $r \leq 
 R_{0},$ the curvature $R$ of $(N, g)$ is uniformly bounded, by (2.15), 
and hence $B_{x_{0}}(R_{0})$ has bounded volume and bounded topology; 
in particular, there is a fixed bound on the number of components of 
$S_{x_{0}}(R_{0}),$ depending only on $R_{0}.$ This proves that the 
number of ends is bounded by the function $\varepsilon (r).$ The same 
reasoning gives an estimate on $V_{0}$ in terms of $\varepsilon (r).$ 

 This completes the proof of Theorem 2.2.
{\endproof}

\begin{remark} \label{r 2.6.}
  {\rm As a simple illustration of Theorem 2.2, let $p_{i}$, $i = 1, 
... \kappa$, be a finite collection of points in a closed Riemannian 
manifold $(M, g)$. Let $r_{i}(x) = dist(x, p_{i})$ and consider the 
manifold $N = M\setminus \cup p_{i}$ with the conformally equivalent 
metric
$$\widetilde g = (\sum_{i=1}^{\kappa}\frac{1}{\widetilde r_{i}^{4}})g, $$
where $\widetilde r_{i}$ is a smoothing of $r_{i}$ away from $p_{i}.$ 
Then $\widetilde g$ is complete and ALE, in fact AE, with $\kappa $ 
ends. The norm of the curvature depends on $\kappa .$ In particular, 
the $L^{n/2}$ and $L^{\infty}$ norms of the curvature are unbounded as 
$\kappa $ increases to infinity. This shows that number of ends, as 
well as the volume ratio $V_{0},$ in Theorem 2.2 depend on $\varepsilon 
(r).$ 

 Similarly, it is easy to construct manifolds satisfying (2.14) and 
(2.15) with $\varepsilon (r) \equiv $ 1 which have infinitely many 
ends. }
\end{remark}

\begin{remark} \label{r 2.7. (i).}
  {\rm It is clear that Theorem 2.2 also holds if $(N, g)$ is an orbifold 
with a bounded number of cones at each singular point, instead of a smooth manifold. 
By definition, the singular points of orbifolds form a discrete set. Hence the cross 
sections $S^{n-1}/\Gamma $ of the tangent cones at infinity cannot themselves be 
orbifolds, for otherwise $N$ would have at least non-trivial curves of singularities. 
Hence, the number of orbifold singular points is finite, and again 
bounded by the function $\varepsilon (r).$

{\bf (ii).}
 In a certain sense, Theorem 2.2 is more general than Theorem 1.2, in 
that the global bound (1.9) on the $L^{n/2}$ norm of the curvature is 
not required in Theorem 2.2. However, for the application to the proof 
of Theorem 1.1, it is important to know that the number $\kappa $ of 
ends, and the maximal volume ratio $V_{0}$ depend only on $\Lambda ,$ 
and not on the manifold $(N, g)$ or the function $\varepsilon (r).$ Of 
course a bound on $\Lambda $ does not give apriori control on 
$\varepsilon (r).$ }
\end{remark}

 Using Theorem 2.2 as a base, we are now in position to prove Theorem 
1.2.

\medskip

{\bf Proof of Theorem 1.2.}

 We recall the hypotheses: $(N, g)$ is a complete non-compact Riemannian 
manifold, satisfying
\begin{equation} \label{e2.51}
vol B_{x}(r) \geq  \nu_{0}r^{n}, \int_{N}|R|^{n/2} \leq  \Lambda , 
\end{equation}
together with the small curvature estimate (1.5), with fixed $\nu_{0}, 
\Lambda , \delta_{0}$ and $C_{0}.$ Let
\begin{equation} \label{e2.52}
{\cal V}_{0}  = \sup_{x, r}\frac{vol B_{x}(r)}{r^{n}}. 
\end{equation}
Then ${\cal V}_{0}$ and $\kappa$, the number of ends of $(N, g)$ are 
finite, and we need to prove there is a bound on ${\cal V}_{0}$ and 
$\kappa$, depending only on $\nu_{0}, \Lambda , \delta_{0}$ and 
$C_{0}.$ Note that all the quantities in (2.51), (2.52) and (1.5) are 
scale-invariant.

 The proof is somewhat simpler if one adds the hypothesis
\begin{equation} \label{e2.53}
|R|(r) \leq  \frac{C}{1+r^{2}}, 
\end{equation}
where $r(x) = dist_{g}(x, x_{0}),$ for some base point $x_{0}\in N$ and 
some arbitrary constant $C < \infty$, with $\kappa$ and ${\cal V}_{0}$ depending 
also on $C$. Thus, we first prove the result in this situation, and then 
show how the proof can be generalized to give the result in full.

 Suppose the statement is false; then there exists a sequence of 
manifolds $(N_{i}, g_{i})$ satisfying (2.51)-(2.53), w.r.t base points 
$x_{i},$ such that either $\kappa_{i} \rightarrow  \infty$ or 
${\cal V}_{i} = {\cal V}_{i}(g_{i}) \rightarrow  \infty$, (or both). 
By the smooth (i.e. $C^{1,\alpha}\cap L^{2,p})$ Gromov compactness theorem, 
there exist radii $R_{i}^{1} \rightarrow  \infty $ such that a 
subsequence of $(B(R_{i}^{1}), g_{i}, x_{i})$ converges uniformly in 
the $C^{1,\alpha}$ and weak $L^{2,p}$ topologies to a complete limit 
$(N_{0}, g_{0}, x_{0}).$ Here $B(R_{i}^{1})$ is the ball centered at 
$x_{i};$ in the following, we assume the subsequence of $\{(N_{i}, 
g_{i})\}$ has been chosen, so that one has convergence to the limit. 
The limit metric $g_{0}$ has curvature locally in $L^{p},$ for any $p < 
 \infty ,$ and has $L^{n/2}$ norm of curvature satisfying (2.51) on 
$N_{0}.$

 We observe that $R_{i}^{1}$ may be chosen so that
\begin{equation} \label{e2.54}
\int_{B(R_{i}^{1})}|R|^{n/2} \geq  \delta_{0}. 
\end{equation}
For if $\int_{B(R_{i}^{1})}|R|^{n/2} <  \delta_{0},$ we may rescale the 
metrics $g_{i}$ by $(R_{i}^{1})^{-2}$ taking $B(R_{i}^{1})$ to a ball 
of radius 1. The small curvature estimate (1.5) implies that the 
curvature is uniformly bounded in $B(1)$, while (2.53) remains valid 
outside $B(1)$; (the estimate (2.53) is essentially scale invariant for 
large $r$). Hence, all the arguments above may be applied to this 
renormalization of the original sequence $(N_{i}, g_{i}).$

 We may thus apply Theorem 2.2 to conclude that $N_{0}$ has a finite 
number $\kappa^{0}$ of ends, and a maximal volume ratio ${\cal V}^{0}.$ 
Returning to the geometry of the original subsequence $(N_{i}, g_{i}),$ 
it follows that a neighborhood $A(\frac{1}{2}R_{i}^{1}, R_{i}^{1})$ of 
the boundary of $B(R_{i}^{1})$ has $\kappa^{0}$ components, and the 
maximal volume ratio ${\cal V}$ in $B(R_{i}^{1})$ satisfies 
${\cal V}  \leq  2{\cal V}^{0}$. The curvature of $g_{i}$ on the region 
$A(\frac{1}{2}R_{i}^{1}, R_{i}^{1}),$ or the larger region 
$A(mR_{i}^{1}, m^{-1}R_{i}^{1}),$ for any fixed $m$, 
is on the order of $\mu_{i}\cdot  (R_{i}^{1})^{-2},$ where 
$\mu_{i} \rightarrow 0$ as $i \rightarrow \infty$; after rescaling 
by $(R_{i}^{1})^{-2},$ the geometry here is close to that of a 
collection of truncated spherical cones.

 Now rescale the original metrics $(N_{i}, g_{i})$ in the subsequence 
to make the radius $R_{i}^{1}$ small. Thus, set $g_{i}^{1} = 
r_{i}^{2}(R_{i}^{1})^{-2}g_{i},$ where $r_{i} \rightarrow $ 0 slowly, 
as $i \rightarrow  \infty .$ The curvature of $g_{i}^{1}$ in the 
$g_{i}^{1}$-annulus $A(mr_{i}, m^{-1}r_{i})$ about $x_{i}$ is on the 
order of $\mu_{i}'  \rightarrow $ 0 as $i \rightarrow  \infty ,$ since 
$r_{i}$ is chosen to go to 0 sufficiently slowly. Outside $B(r_{i}),$ 
the curvature bound (2.53) still holds w.r.t. $g_{i}^{1},$ and $r_{i}$ 
may be chosen so that the curvature is uniformly bounded in 
$(A(r_{i},1), g_{i}^{1}).$ Thus, via Gromov compactness again, by 
passing to a further subsequence of the original sequence, we can again 
pass to a complete rescaled limit $(V^{1}, g^{1}, x^{1})$ with a single 
orbifold singular point $x^{1} = \{0\}$, having a bounded number 
$\kappa^{0}$ of cones. Now apply Theorem 2.2, 
(or Remark 2.7), to this configuration; $V^{1}$ has a bounded number 
$\kappa^{1}$ of ends, and again maximal volume ratio ${\cal V}^{1} 
(\geq  {\cal V}^{0})$. Hence, there is a sequence $R_{i}^{2} \rightarrow  
\infty $ (slowly) such that the balls $B(R_{i}^{2}) \subset  (N_{i}, 
g_{i}^{1})$ converge uniformly to $(V^{1}, g^{1});$ the convergence is 
$C^{1,\alpha}$ and weak $L^{2,p}$ away from the singular point \{0\}. 
For the same reasons as above justifying (2.54), we may assume that, in 
the $g_{i}^{1}$ metric, 
$$\int_{A(1,R_{i}^{2})}|R|^{n/2} \geq  \delta_{0}. $$
Returning then again to the original unscaled subsequence $(N_{i}, 
g_{i}),$ it follows that the number of components of 
$A(\frac{1}{2}R_{i}^{1}R_{i}^{2}, R_{i}^{1}R_{i}^{2})$ is $\kappa^{1}$ 
and the maximal volume ratio ${\cal V} $ in the ball 
$B(R_{i}^{1}R_{i}^{2})$ satisfies ${\cal V}  \leq  {\cal V}^{1}.$

 This process can clearly be repeated any finite number of times, 
passing to a subsequence of the original sequence $(N_{i}, g_{i})$ if 
necessary. Each iteration adds $\delta_{0}$ to the total $L^{n/2}$ norm 
of curvature. Hence, after $d = d(\Lambda , \delta_{0})$ iterations, 
one has a radius $S_{i}^{d} = \prod_{1}^{d}R_{i}^{k}$ such that
\begin{equation} \label{e2.55}
\int_{N_{i}\setminus B(S_{i}^{d})}|R|^{n/2} \leq  \delta_{0}, 
\end{equation}
for which the maximal volume ratio ${\cal V} $ in $B(S_{i}^{d}) \subset 
 (N_{i}, g_{i})$ satisfies ${\cal V}  \leq  {\cal V}^{d} <  \infty $ 
and the number of components of $A(\frac{1}{2}S_{i}^{d}, S_{i}^{d}) 
\subset  (N_{i}, g_{i})$ is bounded by $\kappa^{d} <  \infty ,$ for $i$ 
large. The curvature of the annulus $A(\frac{1}{2}, 2)$ in the 
rescaling $g_{i}^{d} = (S_{i}^{d})^{-2}$ is arbitrarily small when $i$ 
is large. In the scale $g_{i}^{d},$ the estimate (2.55) is now global, 
and implies the curvature decays faster than quadratically outside 
$B(1)$. Theorem 2.2 then gives a uniform bound on the number of ends of 
$N_{i}$ and a uniform upper bound ${\cal V} $ on the maximal volume 
ratio. This contradicts the assumption that ${\cal V}_{i}$ or 
$\kappa_{i}$ becomes arbitrarily large as $i \rightarrow  \infty ,$ and 
proves Theorem 1.2 in case (2.53) holds.

\medskip

 Next we turn to the situation where (2.53) is not assumed. The proof 
in this case has the same overall structure as before; the only 
difference is that new orbifold singular points, not only $\{0\}$, may 
be introduced in the limits. However, these are treated in essentially 
the same way.

 The hypotheses (1.8) and (1.9) in Theorem 1.2 are scale invariant and 
we need to first choose a scale and base points at which to begin the 
analysis. Thus, in the following, we assume that any $(N, g)$ is scaled, 
and a base point $x_{0}\in (N, g)$ is chosen, so that
\begin{equation} \label{e2.56}
\int_{B_{x_{0}}(1)}|R|^{n/2} = \delta_{0}, \ {\rm and} \ 
\int_{B_{x}(1)}|R|^{n/2} \leq  \delta_{0}, \forall x\in N.    
\end{equation}
It follows from the small curvature estimate (1.5) that the curvature 
of $(N, g)$ is uniformly bounded.

 The initial step is then identical to the procedure above, giving a 
limit $(N^{0}, g^{0}, x^{0})$ associated with radii $R_{i}^{1} 
\rightarrow  \infty .$ The first rescaling to obtain the limit $(V^{1}, 
g^{1})$ is also identical within small a small ball $B(\mu_{0})$ about 
the orbifold singularity \{0\}. However, the curvature of $g_{i}^{1}$ 
may not be uniformly bounded outside $B(\mu_{0}).$ Nevertheless, we 
claim that this just corresponds to formation of new orbifold singular 
points on the limit $V^{1}.$ For suppose the curvature of $g_{i}^{1}$ 
blows up at points $y_{i}$ within finite distance to $x_{i}$ w.r.t. 
$g_{i}^{1}.$ By perturbing $y_{i}$ slightly if necessary, we may assume 
that curvature is more concentrated in $L^{n/2}$ norm at $y_{i}$ than 
at nearby points. The metric $g_{i}^{1}$ may then be rescaled to 
$\widetilde g_{i}^{1},$ based at $y_{i}$ so that, with respect to 
$\widetilde g_{i}^{1}$,
$$\int_{B_{y_{i}}(1)}|R|^{n/2} = \delta_{0}, \ {\rm and} \  
\int_{B_{y}(\frac{1}{2})}|R|^{n/2} \leq  \delta_{0},$$
for all $y$ within  $\widetilde g_{i}^{1}$ bounded distance to $y_{i}.$ 
Note that by (2.56), distances w.r.t. the metric  $\widetilde 
g_{i}^{1}$ are no larger than distances w.r.t. the original metric 
$g_{i}.$ One may now repeat the analysis carried out above at $x_{i}$ 
at $y_{i}$ in place of $x_{i}.$ It follows that $y_{i}$ converges, (as 
always in a subsequence of the original sequence) to an orbifold 
singularity $p \in  V^{1}.$ Thus, the limit $(V^{1}, g^{1})$ is a 
complete orbifold, with a bounded number of singular points, (depending 
on $\Lambda ).$ As before, $V^{1}$ has a bounded number $\kappa^{1}$ of 
ends, each ALE, and has a maximal volume ratio ${\cal V}^{1}.$

 With this modification of orbifold limits in place of manifold limits, 
the remainder of the proof is then exactly the same as above in the 
case (2.53) holds. This completes the proof of Theorem 1.2.
{\endproof}

\begin{remark} \label{r 2.8.}
  {\rm As in Remark 2.7, it is clear from the proof above that Theorem 
1.2 also holds for orbifolds with a bounded number of cones at each singular 
point in place of smooth manifolds. It would be of interest to have a direct proof 
of this result, without the use of a contradiction. }
\end{remark}

{\bf \S 2.3.}
 In this section, we complete the proof of Theorem 1.1.

{\bf Proof of Theorem 1.1.}

 Given the work done in the prior subsections, this is now quite 
simple. Thus, let $(M_{i}, g_{i})$ be a sequence of $n$-manifolds 
satisfying (1.3)-(1.5). We claim that 
\begin{equation} \label{e2.57}
vol B_{x}(r) < {\cal V} r^{n}, 
\end{equation}
for some constant ${\cal V}  = {\cal V} (\nu_{0}, \Lambda , \delta_{0}, 
C_{0}),$ where $B_{x}(r)$ is any geodesic $r$-ball in $(M_{i}, g_{i}).$ 
For $r$ sufficiently small, (depending on $i$), one has, for any $x\in 
(M_{i}, g_{i})$, $vol B_{x}(r) \sim  \omega_{n}r^{n}$, since the metric 
is nearly Euclidean on sufficiently small scales. Suppose (2.57) fails, 
so that on some sequence of balls $B_{x_{i}}(r_{i}) \subset  (M_{i}, 
g_{i})$, $vol B_{x_{i}}(r_{i}) \geq  {\cal V}_{i}r_{i}^{n},$ where 
${\cal V}_{i}\rightarrow  \infty $ as $i \rightarrow  \infty .$ Let 
$s_{i}$ be the smallest radius such that, for some $y_{i}\in (M_{i}, 
g_{i}),$
\begin{equation} \label{e2.58}
vol B_{y_{i}}(s_{i}) \geq  2V_{0}s_{i}^{n}, 
\end{equation}
where $V_{0}$ is the $V_{0}$ from Theorem 1.2. If $s_{i} = diam (M, g_{i})$, 
then (2.57) holds, (with ${\cal V}  = 2V_{0}),$ and so there 
is nothing to prove. In fact, if $s_{i}$ is bounded away from 0, 
independent of $i$, then again there is nothing to prove. For in this 
situation, one has 
\begin{equation} \label{e2.59}
vol B_{x}(r) \leq  2V_{0}r^{n}, 
\end{equation}
for all $r \leq  s_{0} \leq  s_{i},$ for some fixed $s_{0} > $ 0. Since 
$vol M_{i} = 1$, this gives an upper bound on the diameter of $(M_{i}, 
g_{i})$ and hence (2.57) holds for all $r \leq diam M$, (with ${\cal 
V}  = CV_{0},$ for a definite $C = C(s_{0})).$ Thus, we may assume that 
$s_{i} \rightarrow $ 0 as $i \rightarrow  \infty .$ Hence, most of 
volume of $(M_{i}, g_{i})$ lies outside $B_{y_{i}}(s_{i}).$ In terms of 
the decomposition (2.3)-(2.4) of $(M_{i}, g_{i}),$ note that one must 
have $y_{i}\in B_{i}^{r},$ for all $i$, and for some $r$ small, tending 
to 0 as $i \rightarrow  \infty .$ 

 Now the scale invariant estimate (2.59) holds on all balls of radius 
$r \leq  s_{i}.$ Rescale the metrics $g_{i}$ by setting $\widetilde 
g_{i} = K^{2}(s_{i})^{-2}g_{i},$ so that the ball $B_{y_{i}}(s_{i})$ 
becomes a ball of radius $K$ w.r.t. $\widetilde g_{i}.$ Here $K$ is a 
fixed large number. Hence, the maximal volume ratio is bounded by 
$2V_{0}$ in all balls of $\widetilde g_{i}$-radius at most $K$. It 
follows from \S 2.1, or the proof of Theorem 1.2 above, that a 
subsequence of $(M_{i}, \widetilde g_{i}, y_{i})$ converges to a 
complete limit orbifold $(V_{\infty}, \widetilde g_{\infty}, 
y_{\infty}),$ satisfying the bounds (1.3)-(1.5). Theorem 1.2 then 
implies that the maximal volume ratio of $(V_{\infty}, \widetilde 
g_{\infty})$ is $V_{0}.$ However, we have $vol B_{y_{\infty}}(K) = 
2V_{0}K^{n}.$ This is a contradiction, which thus proves that (2.57) 
holds.

 Given (2.57), the remainder of the proof of Theorem 1.1 has already 
been completed in \S 2.1.
{\endproof}

\begin{remark} \label{r 2.9.}
  {\rm It is worthwhile to discuss some simple examples illustrating 
Theorem 1.1. As in Remark 2.6, let $p_{j}$, $j = 1, ..., \kappa$, be a 
finite collection of points in a compact Riemannian manifold $(M, g)$ so 
that 
$$\widetilde g = (\sum_{j=1}^{\kappa}\frac{1}{\widetilde r_{j}^{4}})g, $$
is an AE metric on $M \setminus \{p_{j}\}$. Choose a sequence 
$\varepsilon_{i} \rightarrow 0$ as $i \rightarrow  \infty$, and 
consider the metrics 
$$\widetilde g_{i} = \varepsilon_{i}^{2}\widetilde g. $$
As $\varepsilon_{i} \rightarrow 0$, these metrics converge, in the 
Gromov-Hausdorff topology, to a union of $\kappa$ cones $C(S^{n-1})$ 
on $S^{n-1}$, i.e. $\kappa$ balls $B^{n}$, joined at the vertex $\{0\}$. 
The limit metric $g_{\infty}$ is just the flat metric on each ball. The 
boundary of this collection of balls is a collection of $\kappa$ 
spheres $S^{n-1}$, to which one may glue on $\kappa$ balls $B^{n}$ to 
obtain a closed orbifold $V$. Then $V$ is clearly a collection of $\kappa 
$ spheres $S^{n},$ joined at a single vertex. The limit flat metric 
$g_{\infty}$ may be conformally bent, and extended across its boundary 
spheres $S^{n-1}$ to give a smooth metric $g_{\infty}'$ on the 
orbifold $V$. Similarly, the sequence of metrics $\widetilde g_{i}$ may 
easily be modified and extended to give a sequence of smooth metrics 
$\widetilde g_{i}'$ on $M\#(\#_{1}^{\kappa}S^{n}) = M$, converging to 
$g_{\infty}'$ in the Gromov-Hausdorff topology, 
\begin{equation} \label{e2.60}
(M\#(\#_{1}^{\kappa}S^{n}), g_{i}') \rightarrow  
(\vee_{1}^{\kappa}S^{n}, g_{\infty}'). 
\end{equation}
It is easily seen that the sequence $\widetilde g_{i}'$ satisfies all 
the hypotheses of Theorem 1.1. 

 Observe that any non-trivial topology of $M$ has disappeared in the 
limit $V$. It is clear that this construction can be generalized so that 
any closed orbifold $V$ is the limit of a sequence of metrics on 
suitable manifolds $M$ satisfying the hypotheses of Theorem 1.1. The 
only restriction on the topology of $M$ is that the regular part 
$V_{0}$ of $V$ embeds as a domain in $M$, with boundary a union of 
spherical space-forms. The topology of the complement $M\setminus 
V_{0}$ may otherwise be arbitrary, (for $\Lambda $ sufficiently large). 
}
\end{remark}

\begin{remark} \label{r 2.10.}
 {\rm  As outlined in Remark 2.1(ii), one can completely describe the 
structure of the degeneration of $(M_{i}, g_{i})$ near the orbifold 
singularities. We give this description here, but refer to [B] or [AC] 
for more details. The description closely resembles the proof of 
Theorem 1.2.

 Let $(M_{i}, g_{i}) \rightarrow (V, g_{\infty})$ and let $q$ be a 
singular point of $V$. Then there is sequence of scales $r_{i} = 
r_{i}^{1} \rightarrow $ 0 such that the rescalings $(M_{i}, 
r_{i}^{-2}g_{i}, x_{i})$ with $x_{i} \rightarrow q$, converge in the 
pointed Gromov-Hausdorff topology to a complete, non-compact ALE 
orbifold $(V^{1}, g^{1})$ with a finite number of singular points, and
\begin{equation} \label{e2.61}
\int_{V^{1}}|R|^{n/2} \geq  \delta_{0}. 
\end{equation}
If $(V^{1}, g^{1})$ is not a smooth manifold, then there are second 
level scales $\{r_{i}^{2}\}$ associated with each singular point of 
$V^{1};$ (the scales $r_{i}^{2}$ depend on the choice of singular point 
in $V^{1}).$ Rescaling $r_{i}^{-2}g_{i}$ further by such factors at base points 
converging to the singular points of $V^{1}$ gives a collection of second 
level ALE orbifolds $\{(V^{2}, g^{2})\}$ associated with each singular 
point of $(V^{1}, g^{1}).$ Each iteration of this process satisfies 
(2.61), and hence this process terminates at a finite stage. At the last stage, 
corresponding to the smallest scales, the resulting blow-up limits are 
non-flat smooth manifolds. The topology of the original manifolds 
$M_{i}$ may then be reconstructed from that of $V$, and the scale of 
orbifolds associated with each singular point $q\in V$ and its predecessors 
in $V^{j-1}$.

  Each orbifold is of uniformly bounded topological type, with bounds 
depending only on $\Lambda$, $\nu_{o}$, $\delta_{0}$ and $C_{0}$, and 
hence there is a bound $K = K(\Lambda, \nu_{o}, \delta_{0}, C_{0})$ on 
the number of topological types of $\{M_{i}\}$. In particular, 
the homology groups of $M_{i}$ are determined, (by the Mayer-Vietoris 
sequence), by the homology of $V$ and the homology of the collection of 
ALE orbifold spaces $\{V^{j}\}.$ }
\end{remark}

\begin{remark} \label{r 2.11}
{\rm In a certain sense, the bound on $\sup_{B(r/2)}|R|$ in the small curvature 
estimate (1.5) is a stronger assumption than the bound on the Ricci curvature in 
(1.2). For the applications of Theorem 1.1 in \S 3, this turns out to be irrelevant. 
However, it is perhaps worth noting that the bound on $\sup_{B(r/2)}|R|$ can be 
replaced by bound on $\sup_{B(r/2)}|Ric|$, so that (1.5) may be replaced by the 
weaker condition
\begin{equation} \label{e2.62}
\sup_{B(r/2)}|Ric| \leq (\frac{C_{0}}{vol B_{x}(r)}\int_{B_{x}(r)}|R|^{n/2}dV_{g})^{2/n},
\end{equation}
in Theorem 1.1. We outline a proof of this, but do not give all the details, since this 
result is not used.

  To see that (2.62) may be used in place of (1.5), observe that (1.5) is used in the 
proof of Theorem 1.1 in exactly two ways. First, it is used in the local version of the 
$C^{1,\alpha} \cap L^{2,p}$ Gromov compactness theorem. However, as in the proof of 
orbifold compactness under the bounds (1.2), under a smallness condition on the 
$L^{n/2}$ norm of the curvature $R$ and a lower volume bound as in (1.4), such a 
compactness result holds with a bound on $|Ric|$ in place of $|R|$, cf. [An3]. 

  Second, the estimate on $\sup |R|$ is used in the proof of Theorem 1.2 (Steps 1 
and 2) to obtain upper bounds on $\Delta r$ and $D^{2}r$ via comparison geometry. 
These estimates do require a bound on the full curvature $|R|$, not just on the 
Ricci curvature. However, given the bound (1.9) on the $L^{n/2}$ norm of $R$ on 
$(N, g)$, the estimate (2.62) implies that $r^{2}\sup_{A(r,2r)}|Ric| \rightarrow 0$ as 
$r \rightarrow \infty$ on $(N, g)$. Consequently, the $L^{2,p}$ harmonic radius 
$\rho$, (cf. (3.11) below), satisfies $\rho(x) \sim r(x)$ for $r(x)$ large. This 
means that the rescaled metrics $g_{r}$ are locally close to the flat metric in 
the $C^{1,\alpha}$ and weak $L^{2,p}$ topology on balls of uniform size. One may 
then smooth the metric, for instance by convolution with a smooth mollifier, to 
produce a nearby metric for which the estimate (1.5) holds. These local smoothings 
may then be patched together, via a partition of unity, to define a new metric 
$(N, \widetilde g)$ on which (1.5) then holds. Applying Theorem 1.2 to 
$(N, \widetilde g)$ gives then the same conclusion for the original metric $(N, g)$. 

  With these two modifications, the rest of the proof of Theorem 1.1 proceeds in exactly 
the same way. In fact, this discussion shows that (2.62) can be weakened even further, if 
desired. For example, it may be replaced by the bound
\begin{equation} \label{e2.63}
(\frac{1}{vol B_{x}(r/2)}\int_{B_{x}(r/2)}|Ric|^{p}dV_{g})^{1/p} \leq 
(\frac{C_{0}}{vol B_{x}(r)}\int_{B_{x}(r)}|R|^{n/2}dV_{g})^{2/n},
\end{equation}
for a fixed $p > n/2$. One then obtains orbifold compactness in $C^{1,\alpha}$ and 
weak $L^{2,p}$, where $\alpha < 2 - \frac{n}{p}$. }

\end{remark}

\begin{remark} \label{r 2.12} 
{\rm From the proofs of Theorems 1.1 and 1.2, it is apparent that Theorem 1.1 holds for 
spaces of orbifolds satisfying the bounds (1.3)-(1.5). }
\end{remark}

\section{Applications to Moduli of Critical Metrics.}
\setcounter{equation}{0}

 In this section, we prove Theorem 1.3 and related results on metrics 
which are critical points of natural Riemannian functionals on the 
space of metrics. For clarity, we divide this section into two 
subsections. In \S 3.1, we discuss some general results regarding the 
hypotheses of Theorem 1.1. In \S 3.2, we discuss Theorem 1.3, i.e. 
moduli of Bach-flat metrics as well as moduli of selfdual (or 
anti-selfdual) metrics on 4-manifolds, and remark on analogous results 
for critical metrics of other natural functionals on the space of 
metrics.

\medskip

{\bf \S 3.1.}
 This subsection is concerned with some general results relating to the 
hypotheses of Theorem 1.1. 

 Given a closed Riemannian $n$-manifold $(M, g)$, let $c_{S}$ be the 
Sobolev constant for the embedding $L^{2n/n-2} \rightarrow  L^{1,2}$ of 
Sobolev spaces; thus $c_{S}$ is the constant such that
\begin{equation} \label{e3.1}
c_{S}(\int_{M}|f|^{2n/n-2})^{n/n-2} \leq  \int_{M}|df|^{2} + 
\frac{c_{S}}{vol M^{2/n}}\int_{M}f^{2}, 
\end{equation}
for any smooth function $f$ on $(M, g)$. Note that $c_{S}$ is scale-invariant. 
Recall also that a metric $\gamma $ on $M$ is a Yamabe metric if it minimizes the 
total scalar curvature functional restricted to the conformal class $[\gamma]$ of 
$\gamma$; this is equivalent to the statement that, for $n \geq 3$,
\begin{equation} \label{e3.2}
s_{\gamma}vol M^{2/n}(\int_{M}|f|^{2n/n-2})^{n/n-2} \leq  \int_{M}c_{n}|df|^{2} + 
\int_{M}s_{\gamma}f^{2}, 
\end{equation}
for all smooth functions $f$ on $M$, where $c_{n} = 4\frac{n-1}{n-2}$ and 
$s_{\gamma}$ is the scalar curvature of $\gamma$.
\begin{lemma} \label{l 3.1.}
  Let $g$ be a unit volume Yamabe metric on a closed $n$-manifold $M$, with 
scalar curvature $s_{g} \geq s_{0} > 0$. Then there is a constant $\nu_{0}$, 
depending only on $n$ and $s_{0}$, such that
\begin{equation} \label{e3.3}
c_{S} \geq s_{0}/c_{n}, \ \ {\rm and} 
\end{equation}
\begin{equation} \label{e3.4}
vol B_{x}(r) \geq  \nu_{0}r^{n}, 
\end{equation}
for all $r \leq  \frac{1}{2}diam M$.
\end{lemma}
{\bf Proof:}
 Both of these facts are standard and well-known. The estimate (3.3) follows  
trivially from (3.1)-(3.2). The estimate (3.4), while less trivial, is also easy 
to prove, cf. [Ak] for example.
{\endproof}

 Of course the estimate (3.4) corresponds to the non-collapse 
assumption (1.4) in Theorem 1.1. Regarding the global curvature bound 
in (1.4), one has the following, also standard, result.
\begin{lemma} \label{l 3.2.}
  Let (M, g) be a closed (oriented) Riemannian 4-manifold. Then
\begin{equation} \label{e3.5}
\frac{1}{8\pi^{2}}\int_{M}|R|^{2} = \chi (M) + 
\frac{1}{2\pi^{2}}\int_{M}|z|^{2}, 
\end{equation}
\begin{equation} \label{e3.6}
\frac{1}{16\pi^{2}}\int_{M}|z|^{2} = -\chi (M) + 
\frac{1}{8\pi^{2}}\int_{M}|W|^{2} + 
\frac{1}{8\pi^{2}}\int_{M}\frac{s^{2}}{24} , 
\end{equation}
and
\begin{equation} \label{e3.7}
\frac{1}{12\pi^{2}}\int_{M}|W^{+}|^{2} - |W^{-}|^{2} = \tau (M). 
\end{equation}
Here $R$, $W$ and $s$ are the Riemann, Weyl and scalar curvature, $z$ is 
the trace-free Ricci curvature, $z = Ric - \frac{s}{4}g,$ and 
$W^{\pm}$ are the selfdual and anti-selfdual parts of the Weyl 
curvature; $\chi (M)$ and $\tau (M)$ are the Euler characteristic and 
signature.
\end{lemma}
{\bf Proof:}
 These statements (3.5)-(3.6) are just two related versions of the 
Chern-Gauss-Bonnet theorem in dimension 4. The equation (3.7) is a 
combination of the Chern-Weil theorem and the Hirzebruch signature 
formula, cf. [Be].
{\endproof}

 We also recall Aubin's estimate on the scalar curvature of any unit 
volume Yamabe metric $g$ on any $M$:
\begin{equation} \label{e3.8}
s_{g} \leq  S_{0} = n(n-1)\sigma_{n}^{2/n}, 
\end{equation}
where $\sigma_{n}$ is the volume of the unit $n$-sphere in ${\Bbb 
R}^{n+1}.$

 The estimate (3.5) implies an upper bound on the $L^{2}$ norm of the 
curvature $R$ of a metric on a 4-manifold $M$, depending only on an upper 
bound on $\chi (M)$ and an $L^{2}$ bound on $z$; for Einstein metrics, 
one thus obtains an upper bound depending only on $\chi (M).$ 
Similarly, substituting (3.6) in (3.5) and using (3.8) gives an upper 
bound for the $L^{2}$ norm of $R$ on a unit volume Yamabe metric, in 
terms of an upper bound for the $L^{2}$ norm of $W$ and a lower bound 
for $\chi (M).$

 Next we derive an $L^{p}$ analogue of the small curvature estimate 
(1.5).
\begin{lemma} \label{l 3.3.}
  Let B(1) be a geodesic ball centered at $x_{0}$ in a smooth Riemannian 
$n$-manifold $(M, g)$, and suppose the Ricci curvature of $g$ satisfies 
the inequality
\begin{equation} \label{e3.9}
\Delta|Ric| \geq  -d_{0}|R||Ric|, 
\end{equation}
for some constant $d_{0},$ where $R$ is the curvature tensor of (M, g). 
Then there is a constant $\delta_{0},$ depending only on $n$ and the 
Sobolev constant $c_{S}$ in (3.1), such that if 
$(\int_{B(1)}|R|^{n/2})^{2/n} \leq  \delta_{0},$ then
\begin{equation} \label{e3.10}
(\int_{B(\frac{1}{2})}|R|^{p})^{1/p}\leq  C_{0}\delta_{0}, 
\end{equation}
where $C_{0}$ depends only on n, p, $d_{0}$ and $c_{S}.$
\end{lemma}
{\bf Proof:}
 We assume $n$, $p$, $d_{0}$ are chosen, and that $c_{S}$ has a uniform 
lower bound. If $n = 3$, then the Ricci curvature is algebraically 
equivalent to the full curvature $R$, and the proof is much easier in 
this case, cf. Remark 3.4 below. Thus, we assume $n \geq 4$. 

 Let $\rho  = \rho^{2,p}$ be the $L^{2,p}$ harmonic radius of $B(1) 
\subset (M, g)$, cf. [An3]; basically, $\rho $ is the largest radius 
such that on any geodesic ball $B(\rho) \subset B(1)$, one has 
harmonic coordinates in which the metric components $g_{ij}$ are 
bounded by a fixed constant $Q$ in $L^{2,p}$ norm, in that
\begin{equation} \label{e3.11}
Q^{-1}\delta_{ij} \leq  g_{ij} \leq  Q\delta_{ij}, \ \ 
{\rm as \ bilinear \ forms} , 
\end{equation}
$$\rho^{2p-n}\int_{B(\rho )}|\partial^{2}g_{ij}|^{p}dV \leq  Q, $$
where $Q$ is a fixed constant, $Q >$ 1. Thus, $\rho $ is small when the 
curvature is large in $L^{p}$ on $B(\rho )$. The radius $\rho$ scales 
as a distance. Let also $t(x) = dist_{g}(x, \partial B(1))$.

 Consider the minimal value $\mu $ of the (scale invariant) ratio $\rho 
/t.$ The estimate (3.10) clearly holds if $\mu $ is uniformly bounded 
below, i.e. $\rho /t \geq  \mu_{0},$ where $\mu_{0} = 
\mu_{0}(n,p,d_{0}).$ Thus, suppose $\rho (x) << t(x)$ for some $x\in 
B(1),$ so that the $L^{p}$ norm of the curvature becomes large near $x$. 
The minimal value $\mu $ is then achieved at some interior point 
$y_{0}\in B(1),$ and we rescale the metric $g$ at $y_{0}$ to a metric 
$g' $ so that $\rho' (y_{0}) =$ 1, where $\rho' $ is the harmonic 
radius w.r.t. $g' .$ Then $t' (y_{0}) >>$ 1, so that the $g'$-distance 
of $y_{0}$ to the boundary is large. Moreover, the minimality property 
of $y_{0}$ then also gives
\begin{equation} \label{e3.12}
\rho' (y) \geq  \frac{1}{2}, \ \ \forall y\in (B_{y_{0}}(2), g' ). 
\end{equation}

 The equation (3.9) and hypothesis on the $L^{n/2}$ norm of $R$ are 
scale invariant, and so remain valid in the scale $g' .$ For 
simplicity, we drop the prime from the notation, so that one now has 
$\rho (y_{0}) =$ 1, and $\rho (y) > \frac{1}{2},$ for all $y$ within 
distance 2 to $y_{0}.$

 Since the metric is thus uniformly controlled in $L^{2,p} \subset  
C^{1,\alpha}$ on $B(2)$ a standard application of the 
DeGiorgi-Nash-Moser estimate for non-negative solutions of the 
elliptic inequality (3.9), cf. [GT,Thm.8.17], implies that
\begin{equation} \label{e3.13}
sup_{B(\frac{3}{2})}|Ric|^{2} \leq  C\int_{B(2)}|Ric|^{2} + 
C(\int_{B(\frac{3}{2})}|R|^{q})^{2/q}. 
\end{equation}
where $q$ is any fixed number with $q >  n/2$ and $C = C(q, c_{S}).$ 
All balls are now centered at $y_{0}$. A standard interpolation 
inequality, cf. [GT, (7.10)], gives for $n/2 <  q <  p$
\begin{equation} \label{e3.14}
||R||_{L^{q}(B(\frac{3}{2}))} \leq  
\varepsilon||R||_{L^{p}(B(\frac{3}{2}))} + 
\varepsilon^{-\sigma}||R||_{L^{n/2}(B(\frac{3}{2}))}, 
\end{equation}
for any $\varepsilon  > $ 0, where $\sigma  = (2/n$ - $1/q)/(1/q$ - 
$1/p).$  Also, from the fact that $\rho  \geq  \frac{1}{2}$ in 
$B(2)$, we have
\begin{equation} \label{e3.15}
||R||_{L^{p}(B(\frac{3}{2}))} \leq  C||R||_{L^{p}(B(1))},  
\end{equation}
for a fixed constant $C$.

 On the other hand, the Ricci curvature controls the metric in harmonic 
coordinates, and one also has the elliptic estimate
\begin{equation} \label{e3.16}
||R||_{L^{p}(B(1))} \leq  C_{1}(||Ric||_{L^{\infty}(B(\frac{3}{2}))} + 
||R||_{L^{n/2}(B(\frac{3}{2}))}). 
\end{equation}
coming from the well-known elliptic equation $(Ric)_{ij} = 
-\frac{1}{2}\Delta g_{ij} + Q(g, \partial g),$ (cf. [An3,(2.7)ff]). 
Choosing $\varepsilon$ small in (3.14), the combined equations 
(3.13)-(3.16) give the estimate
\begin{equation} \label{e3.17}
||R||_{L^{p}(B(1))} \leq  C_{2}(||Ric||_{L^{2}(B(2))} + 
||R||_{L^{n/2}(B(\frac{3}{2}))}), 
\end{equation}
for some fixed constant $C_{2}.$ Now the last term in (3.17) is bounded 
by $C_{2}\cdot \delta_{0}.$ Similarly, since $|Ric| \leq  |R|,$ and 
since $n \geq $ 4, the $L^{2}$ norm of the Ricci curvature is bounded 
by the $L^{n/2}$ norm of the full curvature $R$. Thus (3.17) implies
\begin{equation} \label{e3.18}
||R||_{L^{p}(B(1))} \leq  C_{3}\delta_{0}. 
\end{equation}
However, the $L^{p}$ norm of $R$ on $B(1)$ is on the order of $Q$ from 
(3.11). Hence, one has a contradiction if $\delta_{0}$ is sufficiently 
small.
{\endproof}

\begin{remark} \label{r 3.4.(i).}
  {\rm If instead of (3.9) one has the inequality $\Delta|R| \geq  
-d_{0}|R|^{2},$ then (3.10) also holds. The proof in this case is a 
direct consequence of the DeGiorgi-Nash-Moser theorem, cf. [An1], 
[BKN]. Such an estimate holds for the curvature tensor of Einstein 
metrics, or more generally metrics with harmonic curvature, cf. [Be]. 

{\bf (ii).}
 The same proof shows that (3.10) if the hypothesis on the $L^{n/2}$ 
norm is replaced by other scale invariant bounds on the curvature. For 
example, the proof also holds in case one assumes
\begin{equation} \label{e3.19}
\sup_{B(r)\subset B(1)}\frac{1}{r^{n-2s}}\int_{B(r)}|R|^{s} \leq  
\delta_{0}, 
\end{equation}
for 2 $\leq  s \leq  n/2.$ }
\end{remark}

 Lemma 3.3 leads to the following strong or smooth version of the small 
curvature estimate (1.5).
\begin{proposition} \label{p 3.5.}
  Let B(1) be a geodesic ball in a smooth Riemannian n-manifold (M, g) 
and suppose the Ricci curvature of $g$ satisfies an equation of the form
\begin{equation} \label{e3.20}
D^{*}DRic = R*Ric, 
\end{equation}
where $R$ is the curvature tensor and $*$ denotes any algebraic 
operation of $R$ or its (irreducible) components on symmetric bilinear 
forms. Suppose in addition that $(\int_{B(1)}|R|^{n/2})^{2/n} \leq  
\delta_{0},$ where $\delta_{0}$ is small, depending only on $n$ and 
$c_{S}.$

 Then there is an $r_{0} = r_{0}(n, c_{S})$ and harmonic coordinates on 
$B(r_{0})$ such that
\begin{equation} \label{e3.21}
||g||_{C^{m,\alpha}(B(r_{0}))} \leq  K, 
\end{equation}
where $K$ depends only on n, m, $\alpha ,$ and $c_{S}.$ 
\end{proposition}
{\bf Proof:}
 Taking the inner product of (3.20) with $Ric$ gives 
$$\tfrac{1}{2}\Delta|Ric|^{2} = |DRic|^{2} + \langle R*Ric, Ric \rangle.$$
A simple application of the Cauchy-Schwartz inequality then implies
$$\Delta|Ric| \geq  -d_{0}|R||Ric|, $$
for a constant $d_{0}$ depending only on the algebraic operator $*.$ 
From Lemma 3.3, one thus has a uniform lower bound on the $L^{2,p}$ and 
therefore $C^{1,\alpha}$ harmonic radius. In particular, (3.21) holds 
with $m =$ 1. 

 Next, since the curvature is uniformly bounded in $L^{p}$ on 
$B(\frac{1}{2}),$ standard elliptic regularity, [GT, Ch.8] applied to the 
operator (3.20) expressed in local harmonic coordinates on $B(r_{0}),$ 
(in which the metric is $L^{2,p})$ implies that $Ric$ is bounded in 
$L^{4,p/2}$ on $B(r_{0}/2)$. This, together with the fact that the 
$L^{2,p}$ harmonic radius is bounded below by $r_{0}$ implies that the 
$L^{4,p/2}$ harmonic radius is bounded below by $r_{0}/2$. It follows 
that the metric is controlled in $L^{4,p/2}$, and so in $C^{3,\alpha}$ on 
$B(r_{0}/2)$. By taking a suitable covering, it follows that the metric 
is controlled in $C^{3,\alpha}$ on $B(r_{0}).$ This process is then 
continued inductively to obtain (3.21).
{\endproof}

 Finally, we conclude with a remark on the Sobolev inequality. The 
Sobolev inequality implies a lower bound on the volume ratio of 
geodesic balls, as in (3.4). As is well-known, the converse is not 
true in general. This has the following consequence:
\begin{proposition} \label{p 3.6.}
  Suppose $(M_{i}, g_{i})$ satisfy hypotheses of Theorem 1.1, with 
(1.4) replaced by a lower bound on the Sobolev constant $c_{S}.$ Let 
$V$ be a limit orbifold of $(M_{i}, g_{i}),$ so that any spherical 
space form $S^{n-1}/\Gamma_{j}$ in a cone $C(S^{n-1}/\Gamma_{j})$ near 
a singular point $q\in V$ embeds in $M_{i},$ for $i$ large. 

 If some space-form $S^{n-1}/\Gamma $ near $q$ bounds a domain $U_{i} 
\subset  M_{i},$ then the orbifold singularity $q$ is irreducible, so 
that $C(S^{n-1}/\Gamma )$ is the only cone at $q$ in V.
\end{proposition}
{\bf Proof:}
The proof is standard, although we include a proof for completeness. 
In the following, we assume that $i$ is sufficiently large, but 
drop the index $i$ from the notation. Since $(M, g)$ is close to $V$, and 
$S^{n-1}/\Gamma$ separates $M$, one has a decomposition $M = U_{1} \cup U_{2}$, 
where the union is along a neck region $A(\varepsilon, 2\varepsilon)$ whose 
geometry, when scaled by $\varepsilon^{-1}$ is close to that of a standard 
annulus $A(1,2)$ in a standard cone $C(S^{n-1}/\Gamma)$. In particular, 
$vol A(\varepsilon, 2\varepsilon) \sim \varepsilon^{n}$. Let 
$v_{1} = vol U_{1}$, $v_{2} = vol U_{2}$ and assume $v_{1} \leq v_{2}$, 
so that $v_{1} \leq \frac{1}{2}$. Let $f$ be a function on $M$ with 
$f = 1$ on $U_{1} \setminus A(\varepsilon, 2\varepsilon)$, 
$f = 0$ on $U_{2} \setminus A(\varepsilon, 2\varepsilon)$, with 
$|df| \sim \varepsilon^{-1}$ on $A(\varepsilon, 2\varepsilon)$. 
Using $f$ as a test function in the Sobolev inequality (3.1) gives
$$c_{S}v_{1}^{(n-2)/n} \leq c\varepsilon^{n-2} + c_{S}v_{1}.$$
Since $v_{1} \leq \frac{1}{2}$, $v_{1} \leq (1/2)^{2/n}v_{1}^{(n-2)/n}$, 
and hence one has
$$c_{S}v_{1}^{(n-2)/n}(1 - (\tfrac{1}{2})^{2/n}) \leq c\varepsilon^{n-2}.$$
Since $\varepsilon$ may be chosen to be arbitrarily small, (for $i$ sufficiently 
large), $v_{1}$ is arbitrarily small, for $i$ sufficiently large.  It follows easily 
that $V$ contains only a single cone at $q$.
{\endproof}

 The hypothesis that $S^{n-1}/\Gamma $ bounds a domain in $M$ (or 
$M_{i})$ is automatically the case if for instance the Betti number 
$b_{n-1}(M) =$ 0, or equivalently $b_{1}(M) = 0$, by Poincar\'e duality. 

\begin{remark} \label{r 3.7.}
{\rm  A lower bound on the Ricci curvature bound and upper and lower 
bounds on the volume of $(M, g)$ imply a lower bound on the Sobolev 
constant, cf. [Cr], [Ga]. Hence, these bounds imply irreducibility of 
all orbifold limits in Theorem 1.1. On the other hand, Remark 2.9 shows 
that in general, one may have many cones joined at a single orbifold 
singularity.

 An interesting example illlustrating the fact that $S^{n-1}/\Gamma $ 
must bound a domain in Proposition 3.6 is the following. Let $M = M_{i} 
= S^{n-1}\times S^{1}.$ There exists a sequence of conformally flat, unit 
volume Yamabe metrics $g_{i}$ on $M$ which satisfy the hypotheses of 
Theorem 1.1 and whose scalar curvatures tend to the maximal value 
$S_{0}$ in (3.8), cf. [Ko], [Sc]. Hence, there is a uniform Sobolev 
inequality on $(M, g_{i})$ The limit $(V, g_{\infty})$ is the round 
unit volume metric on $S^{n},$ with two antipodal points 
identified $\{q\} \sim \{-q\}$; at this point, one then has an orbifold 
singularity, consisting of two cones, each a ball, identified at the 
origin. }
\end{remark}

{\bf \S 3.2.}
 We now apply the results above to critical metrics for natural 
functionals on the space of metrics. 

{\bf Proof of Theorem 1.3.}

 Let ${\cal C}^{+}$ be the space of unit volume Bach-flat Yamabe 
metrics on a given 4-manifold $M$, with scalar curvature $s_{g} \geq  
s_{0}.$ Theorem 1.3 follows from Theorem 1.1 if we show that the 
hypotheses of the latter are satisfied. 

 First, the non-collapse estimate (1.4) follows from (1.14) and (3.4). 
The global curvature bound (1.3) follows from the bound (1.15) on the 
Weyl curvature, together the expression
\begin{equation} \label{e3.22}
\frac{1}{8\pi^{2}}\int_{M}|R|^{2} = -7\chi (M) + 
\frac{1}{\pi^{2}}\int_{M}|W|^{2} + 
\frac{1}{\pi^{2}}\int_{M}\frac{s^{2}}{24} , 
\end{equation}
valid on any closed 4-manifold $(M, g)$, obtained by combining 
(3.5)-(3.6). The volume is normalized to 1, and so the scalar curvature 
satisfies (3.8).

 To obtain the small curvature estimate, observe that Bach-flat Yamabe 
metrics satisfy an equation of the form (3.20). Namely, the Bach 
equation (1.13) together with the Weitzenbock formula on vector valued 
2-forms gives the equivalent equation
\begin{equation} \label{e3.23}
D^{*}DRic + \tfrac{1}{3}D^{2}s + \tfrac{1}{6}\Delta s\cdot  g + {\cal R}  
= 0, 
\end{equation}
where ${\cal R} = R*Ric$. This equation is conformally invariant. However, 
in the Yamabe normalization, $s = const$, so the Hessian and Laplace terms 
vanish, so that
\begin{equation} \label{e3.24}
D^{*}DRic + {\cal R}  = 0, 
\end{equation}
i.e. (3.20) holds for Bach-flat Yamabe metrics. The small curvature 
estimate (1.5) is then a consequence of Proposition 3.5. Theorem 1.3 is 
thus a consequence of Theorem 1.1.
{\endproof}

\begin{remark} \label{r 3.8.}
 {\rm  With some further work, it can probably be shown that a 
Bach-flat Yamabe metric extends smoothly $(C^{\infty})$ over \{0\} in a local 
uniformization $B^{n}\setminus \{0\}$ of $C(S^{n-1}/\Gamma) \setminus 
\{q\}$. By means of a standard cutoff function argument, it is 
straightforward to see that $g$ is a weak solution of the equation 
(3.21) on the full ball $B^{n}.$ }
\end{remark}

 It is clear that Theorem 1.3 also holds if one allows the manifold $M$ 
to vary, provided, via (3.22), one assumes a uniform lower bound on the 
Euler characteristic of $M$. Again, we point out that the $L^{2}$ norm 
${\cal W} $ of $W$ is constant on components of moduli spaces of 
Bach-flat metrics. Also, by Proposition 3.6, the limit orbifold is 
necessarily irreducible if $b_{1}(M) =$ 0, for instance if $M$ is 
simply connected. (Needless to say, Theorem 1.3 also holds for 
Bach-flat Yamabe metrics of non-positive scalar curvature, provided one 
assumes the bound (1.4)).

\medskip

 Examples of Bach-flat metrics, besides Einstein metrics, (which are 
necessarily Yamabe by Obata's theorem, cf. [Be,Ch.4]), are conformally 
flat metrics, and half-conformally flat metrics. A considerable theory 
for the existence of half-conformally flat metrics has been developed 
by LeBrun, Taubes, and others, cf. [Le], [Ta] and references therein. 
Theorem 1.3 may be used to study aspects of the moduli spaces of such 
metrics. 

 For conformally flat metrics, (3.22) gives an upper bound on $L^{2}$ 
norm of $R$, depending only on a lower bound for $\chi (M).$ Hence, one 
has orbifold compactness, just depending on a lower bound for $\chi 
(M).$ A similar result holds for selfdual (or anti-selfdual) metrics. 
Thus suppose $W = W^{+},$ so that $W^{-} = 0$. Substituting (3.7) in 
(3.22) then gives, for selfdual metrics,  
\begin{equation} \label{e3.25}
\frac{1}{8\pi^{2}}\int_{M}|R|^{2} = -7\chi (M) + 12\tau (M) + 
\frac{1}{\pi^{2}}\int_{M}\frac{s^{2}}{24} . 
\end{equation}
Thus, one again has an upper bound on the $L^{2}$ norm of curvature, 
depending only on an upper bound for $-7\chi (M)+12\tau (M).$

 The formulas (3.5)-(3.7) have analogues for complete, non-compact ALE 
manifolds, (or orbifolds) $(N, g)$ satisfying the volume growth condition 
(1.8). All such ALE spaces which arise as blow-up limits of 
degenerating Bach-flat Yamabe metrics are necessarily scalar-flat. For 
simplicity, we thus state the formulas for scalar-flat ALE manifolds:
\begin{equation} \label{e3.26}
\frac{1}{8\pi^{2}}\int_{N}|R|^{2} = -7\chi (N) + 
\frac{1}{\pi^{2}}\int_{N}|W|^{2} + 
7\sum_{k=1}^{\kappa}\frac{1}{|\Gamma_{k}|}, 
\end{equation}

\begin{equation} \label{e3.27}
\frac{1}{12\pi^{2}}\int_{N}|W^{+}|^{2} - |W^{-}|^{2} = \tau (M) + 
\sum_{k=1}^{\kappa}\eta (S^{3}/\Gamma_{k}), 
\end{equation}
and, in the selfdual case
\begin{equation} \label{e3.28}
\frac{1}{8\pi^{2}}\int_{N}|R|^{2} = -7\chi (N) + 12\tau (N) + 
7\sum_{k=1}^{\kappa}\frac{1}{|\Gamma_{k}|} + 12\sum_{k=1}^{\kappa}\eta 
(S^{3}/\Gamma_{k}). 
\end{equation}

\begin{corollary} \label{c 3.9.}
  Let $g$ be a complete self-dual, scalar-flat, metric on ${\Bbb R}^{4}$ 
satisfying
\begin{equation} \label{e3.29}
\int_{N}|R|^{2} <  \infty , \ \ vol B_{x}(r) \geq  \nu_{0}r^{4}. 
\end{equation}
Then $g$ is flat.
\end{corollary}
{\bf Proof:}
 The metric $({\Bbb R}^{4},$ g) satisfies all the hypotheses of Theorem 
1.2, and hence (3.28) holds. Since in this case $\kappa  = 1$, $\chi 
({\Bbb R}^{4}) =$ 1, $\tau ({\Bbb R}^{4}) = 0$ and $\Gamma = \{e\}$, 
the result follows. An alternate proof is to use (3.27) to see that 
$g$ is conformally flat. The Liouville theorem implies that a scalar-flat 
conformally flat metric on ${\Bbb R}^{4}$ is flat.
{\endproof}

 This result can be viewed as an analogue of the well-known result, 
(which follows from (3.7)) that the only self-dual metric on 
$S^{4}$ is the conformal class of the round metric. The condition 
(3.29) on the volume growth is necessary. For example, the Taub-NUT 
metric [Be,Ch.13] is a complete finite action selfdual Ricci-flat metric 
on ${\Bbb R}^{4}$ which is not flat.
\begin{remark} \label{r 3.10.}
 {\rm  Versions of Theorem 1.3 also hold for other natural functionals 
on the space of metrics. For example, in dimension 4, the 
Euler-Lagrange equations for the $L^{2}$ norms of the full curvature $R$, 
the Ricci curvature $Ric$ or the trace-free Ricci curvature $z$,
\begin{equation} \label{e3.30}
{\cal R}^{2} = \int_{M}|R|^{2}dV,  \ \ {\cal R}ic^{2} = 
\int_{M}|Ric|^{2}dV, \ \  {\cal Z}^{2} = \int_{M}|z|^{2}dV, 
\end{equation}
are all of the form (3.20), cf. [Be,Ch.4]. Hence, any metric which is a 
critical point of such a functional satisfies the small curvature 
estimate (1.5). Such critical metrics also have constant scalar 
curvature; however, it is unknown if they are Yamabe metrics. Thus, it 
is not known if the lower volume bound (1.4) follows from a positive 
lower bound on the scalar curvature, and so (1.4) must be assumed. The 
bound (1.3) on the $L^{2}$ norm of the curvature follows from a bound 
on any of these functionals via (3.5). Under these hypotheses then, 
Theorem 1.3 also holds for moduli spaces of critical metrics of the 
functionals (3.30).

 It is an interesting open problem whether Theorem 1.3 can be extended 
toward an understanding of the degeneration of Bach-flat Yamabe metrics 
with non-positive scalar curvature. This will necessarily include 
collapsing families of metrics, in the sense of Cheeger-Gromov, and the 
formation of limits with cusp-like ends. For Einstein metrics, such a 
study was carried out in [An4]. }
\end{remark}

\bibliographystyle{plain}

\bigskip
\begin{center}
December, 2003/January, 2004
\end{center}

\medskip
\noindent
\address{Department of Mathematics\\
SUNY at Stony Brook\\
Stony Brook, NY 11794-3651}\\
\email{anderson@@math.sunysb.edu}

\end{document}